\documentclass[11pt]{amsart}
\usepackage{amssymb, adjustbox, enumerate, amsbsy, stmaryrd}
\usepackage{amsthm}
\usepackage{geometry}
\usepackage[all]{xy}
\usepackage{ulem}
\usepackage{todonotes}

\usepackage{amsfonts, amssymb, amscd}
\numberwithin{equation}{section}

\usepackage{bm}
\usepackage{verbatim}
\usepackage{mathrsfs}
\usepackage{graphicx}
\usepackage{tikz-cd}
\usepackage{subcaption}
\usepackage{listings}
\usepackage{subfiles}
\usepackage[toc,page]{appendix}
\usepackage{mathtools}
\usepackage{comment}
\usepackage{enumerate}
\usepackage{enumitem}
\usepackage[all]{xy}

\usepackage{graphicx}
\graphicspath{{images/}}

\usepackage{appendix}
\usepackage{hyperref}
\hypersetup{
    colorlinks=true,
    citecolor=red,
    linkcolor=blue,
    filecolor=magenta,      
    urlcolor=red,
}
\def\ideal#1.{I_{#1}}
\def\ring#1.{\mathcal {O}_{#1}}
\def\fring#1.{\hat{\mathcal {O}}_{#1}}
\def\proj#1.{\mathbb P(#1)}
\def\pr #1.{\mathbb P^{#1}}
\def\af #1.{\mathbb A^{#1}}
\def\Hz #1.{\mathbb F_{#1}}
\def\Hbz #1.{\overline{\mathbb F}_{#1}}
\def\pic#1.{\operatorname {Pic}\,(#1)}
\def\pico#1.{\operatorname{Pic}^0(#1)}
\def\picg#1.{\operatorname {Pic}^G(#1)}
\def\ner#1.{NS (#1)}
\def\rdown#1.{\llcorner#1\lrcorner}
\def\rup#1.{\ulcorner#1\urcorner}
\def\cone#1.{\operatorname {NE}(#1)}

\def\ccone#1.{\overline{\operatorname {NE}}(#1)}
\def\coef#1.{\frac{(#1-1)}{#1}}
\def\vit#1.{D_{\langle #1 \rangle}}
\def\mm#1.{\overline {M}_{0,#1}}
\def\H1#1.{H^1(#1,{\ring #1.})}
\def\ac#1.{\overline {\mathbb F}_{#1}}

\def\adj#1.{\frac {#1-1}{#1}}
\def\spn#1.{\overline{#1}}
\def\ses#1.#2.#3.{0\to #1\to #2\to #3 \to 0}
\def\pek#1.#2.{\Cal P^{#1}(#2)}
\def\plk#1.#2.{\Cal P^{\leq #1}(#2)}
\def\ev#1.{\operatorname{ev_{#1}}}
\def\bminv#1.{(\nu_1,s_1;\nu_2,s_2;\dots ;\nu_{#1},s_{#1};\nu_{r+1})}
\def\zinv#1.{(\nu_1,s_1;\nu_2,s_2;\dots ;\nu_{#1},s_{#1};0)}
\def\iinv#1.{(\nu_1,s_1;\nu_2,s_2;\dots ;\nu_{#1},s_{#1};\infty)}
\def\map#1.#2.{#1 \longrightarrow #2}
\def\rmap#1.#2.{#1 \dasharrow #2}
\def\emb#1.#2.{#1 \hookrightarrow #2}

\def\dim{\operatorname{dim}}

\def\bbeta{\boldsymbol{\beta}}

\def\C{\mathbb C}

\def\Z{\mathbb Z}

\def\e{\Cal E}

\def\e1{E_1}
\def\e2{E_2}

\def\OO{\mathcal O}
\def\LL{\mathcal L}

\newcommand{\po}{\ar@{}[dr]|{\text{\pigpenfont R}}}
\newcommand{\pb}{\ar@{}[dr]|{\text{\pigpenfont J}}}

\newcommand\Q{{\mathbb{Q}}}
\newcommand\R{{\mathbb{R}}}

\newtheorem{thm}{Theorem}[section]

\newtheorem{claim}[thm]{Claim}

\theoremstyle{definition}

\theoremstyle{definition}

\newtheorem{theorem}{Theorem}[section]
\newtheorem{lemma}[theorem]{Lemma}
\newtheorem{proposition}[theorem]{Proposition}
\newtheorem{corollary}[theorem]{Corollary}

\theoremstyle{definition}

\newtheorem{conjecture}[theorem]{Conjecture}

\begin{document}
\author{Christopher Hacon} \address{Department of Mathematics\\ University of Utah\\ 155 S 1400 E\\ Salt Lake City, Utah 84112} \email{hacon@math.utah.edu}
\thanks{Christopher Hacon was partially supported by the NSF research grants no: DMS-1952522, DMS-1840190, DMS-2301374. }
\author{Pengchao Wang} 
\address{Academy of Mathematics and Systems Science, Chinese Academy of Sciences,}
\email{blowa8@gmail.com}
\date{\today}
\title[Termination of the K\"ahler MMP with scaling]{Termination of the K\"ahler MMP with scaling}

\begin{abstract}
    In this article,   we prove termination of klt flips with scaling for strongly \(\mathbb Q\)-factorial compact Kähler klt pairs $(X,B)$ with effective \(\mathbb Q\)-boundary. In particular, compact Kähler klt pairs with effective \(\mathbb Q\)-boundary have minimal models when \(K_X+B\) is pseudo-effective, and admit Mori fibre space models otherwise.\end{abstract}
    \maketitle
\tableofcontents
\section{Introduction} 
In recent years there has been dramatic progress towards the minimal model program for K\"ahler varieties.
In \cite{HX26} the authors prove the existence of flips, flipping and divisorial contractions and the termination of flips with scaling and  existence of minimal models for pseudo-effective generalized klt pairs $(X,B+\bbeta )$ where $B+\bbeta$ is modified big and the existence of Mori fiber spaces for non-pseudo-effective generalized klt pairs.  In this paper we prove the termination of flips with scaling for klt pairs and the existence of minimal models (or Mori fiber spaces) for the usual pairs, without any bigness assumption.
\begin{theorem}\label{t-main} Let $(X,B)$ be a strongly $\Q$-factorial compact K\"ahler klt pair, $B$ an effective $\Q$-divisor, and $\omega$ a K\"ahler class such that $K_X+B+\omega$ is nef. Then any $K_X+B$ minimal model program with scaling of $\omega$ terminates after finitely many steps.
\end{theorem}
Combining Theorem \ref{t-main} and strong \(\mathbb Q\)-factorialization with \cite[Theorem 1.2]{OAI3}, where subadditivity assumption is supplied by \cite[Theorem 1.2]{OAI4}, gives the following consequence.
\begin{corollary}\label{c-1.2} Let $(X,B)$ be a compact K\"ahler klt pair, where \(B\) is an effective \(\mathbb Q\)-divisor. If $K_X+B$ is pseudo-effective, then $(X,B)$ has a good minimal model, in particular $\kappa (K_X+B)\geq 0$.
\end{corollary}
We remark that the termination of flips with scaling for generalized klt pairs is still an open conjecture. More generally one expects the following.
\begin{conjecture} Let $(X,B+\bbeta)$ be a strongly $\Q$-factorial compact K\"ahler generalized log canonical pair. Then there are no infinite sequences of $K_X+B+\bbeta _X$ flips. 
\end{conjecture}
We do not discuss the history of the existence of minimal models and the termination of flips, however we refer the reader to the following papers and references therein  \cite{BCHM10}, \cite{HQZ25}, \cite{Has26}, \cite{OAI1}, \cite{OAI2}.

\subsection{Sketch of the proof}\label{ss-sketch}
If \(K_X+B\) is not pseudo-effective, then \cite[Theorem 5.16(2)]{HX26} implies that the given \((K_X+B)\)-MMP with scaling of \(\omega\) terminates with a Mori fibre space. We therefore assume that \(K_X+B\) is pseudo-effective.
There are two main ingredients in the proof of Theorem \ref{t-main} corresponding to a bound on the discrepancies and on the Cartier indices. In order to explain the proof we fix some notation. Let $X,B,\omega $ be as in Theorem \ref{t-main}. By \cite{HX26}, we may run the mmp with scaling of the K\"ahler form $\omega$. More precisely there exists numbers  $1=\lambda _0\geq \lambda _1\geq \lambda _2\geq \ldots \geq 0$, and a sequence of $K_X+B$ flips and divisorial contractions \[X\dasharrow X_1\dasharrow  X_2 \dasharrow \ldots \] such that if $\phi _i:X\dasharrow X_i$, $B_i=(\phi _i)_*B$, and $\omega _i=(\phi _i)_*\omega$, then $X_i$ is strongly $\Q$-factorial, $K_{X_i}+B_i+t\omega _i$ is nef for $\lambda _i\geq t\geq \lambda _{i+1}$.
For convenience we set $X=X_0$, $B=B_0$, $\omega =\omega _0$, and $\lambda _0=1$.
If $\psi _i:X_i\dasharrow X_{i+1}$ is a flip, then we let $f_i:X_i\to Z_i$ be the flipping contraction and $f_i^+:X_{i+1}\to Z_i$ be the flipped contraction. We aim to show that there is no infinite sequence of such maps with positive thresholds $\lambda _i>0$. Arguing by contradiction, we will suppose to the contrary that this sequence is infinite. 

We observe that $\lambda:=\lim \lambda _i=0$. If this is not the case, 
then let \(\lambda_\infty:=\lim\lambda_i>0\), and choose \(0<\delta<\lambda_\infty\). The given sequence is then a \((K_X+B+\delta\omega)\)-MMP with scaling of \((1-\delta)\omega\), 
and scaling numbers \((\lambda_i-\delta)/(1-\delta)\).  Viewing \(\delta\omega\) as a nef \(b\)-class descending to \(X\), the induced generalized pair $(X,B+\delta \bar \omega)$ is gklt and its generalized boundary is big. This contradicts \cite[Theorem 1.5]{HX26}, and therefore \(\lambda_\infty=0\). 

It is easy to see that this sequence contains only finitely many divisorial contractions, and so we may assume that for $i\gg 0$, each $\psi _i$ is a flip.
One of the key ingredients in the proof is the following.
\begin{proposition}\label{p-1} Let $X,B,\omega$ be as in Theorem \ref{t-main}. If $X\dasharrow X_1\dasharrow  X_2 \dasharrow \ldots $ is a sequence of steps of the $K_X+B$ mmp with scaling of $\omega$, then there exists an integer $R>0$ such that $R(K_{X_i}+B_i)$ is Cartier for any $i>0$.
\end{proposition}
Since $-(K_{X_i}+B_i)$ is ample over $Z_i$, then by Koll\'ar’s effective basepoint-freeness (see \cite[Theorem 6.4]{Fuj22}) there is an integer $N >0$ depending
only on $\dim X$ such that $-NR(K_{X_i}+B_i)$ is $f_i$-generated and
$f_i$-ample. 
Suppose that for $i\gg 0$, we have \begin{equation}\label{e-rv}R^j{f_i}_*\OO _{X_i}(m_j(K_{X_i}+B_i))=0\qquad {\rm for}\ 0<j\leq n,\ m_j=R+RNj,  \end{equation} then by Mumford regularity $\OO _{X_i}(R(K_{X_i}+B_i))$ is $f_i$-generated (see \cite[Example 1.8.24]{Laz04}), contradicting the fact that $-(K_{X_i}+B_i)$ is $f_i$-ample.

In order to prove the relative vanishing, note that since $K_{X_{i+1}}+B_{i+1}$ is $f_i^+$-ample and klt, then by relative Kawamata-Viehweg vanishing, $R^j{f_i^+}_*\OO _{X_{i+1}}(m_j(K_{X_{i+1}}+B_{i+1}))=0$.
The above vanishing \eqref{e-rv}, then follows from a delicate discrepancy bound (see Proposition \ref{p-5} and equation  \eqref{e-8}) that allows us to compare the pullbacks of $\OO _{X_i}(m_j(K_{X_i}+B_i))$ and $\OO _{X_{i+1}}(m_j(K_{X_{i+1}}+B_{i+1}))$ on a common resolution.

The proof of Proposition \ref{p-1} is quite delicate. 
 Since the claim is obvious for a finite mmp, we will assume that our mmp does not terminate.
Fix \(p:W\to X_0\) a log resolution of $(X,B)$, a positive integer $d$ such that $dB$ is integral, and put $D_i=(K_{X_i}+B_i)$, $L_i:=dD_i$, and $A_i(v):=A_{X_i,B_i}(v)$.
It suffices to uniformly bound the local Cartier indices of the integral $\mathbb Q$-Cartier divisors $L_i$. Let $E_i:=A_i-A_0$, and choose a positive integer $k$ such that $kdD_0$ is Cartier on $X=X_0$.

Suppose, by way of contradiction, that the Cartier indices are unbounded.
We can choose a subsequence $i_j\to\infty$ and points $x_j\in X_{i_j}$ such that
\[r_j:=\operatorname{ind}_{x_j}(L_{i_j})\longrightarrow\infty. \]
By Lemma \ref{lem:gorenstein-index-cover}, there are finite quasi-\'etale Galois morphisms of pointed germs
\[\pi_j:(V_j,o_j)\longrightarrow(X_{i_j},x_j)\]
with groups $H_j$ fixing $o_j$, such that $V_j$ is Gorenstein klt and $\pi_j^*L_{i_j}$ is Cartier. Moreover, $\mathcal O_{V_j}(\pi_j^*L_{i_j})$ has a semi-invariant local generator $\eta_j$ whose character $\chi_j:H_j\to\mu_{r_j}$ is surjective.
Fix an arbitrary $\alpha\in[0,1)$, set
\[a_j:=\frac{\lfloor r_j\alpha\rfloor}{r_j},\] and choose $g_j\in H_j$ with $\chi_j(g_j)=\exp(2\pi {\rm i} a_j)$. We then have that $\lim _{j\to \infty} a_j=\alpha$.

Applying Proposition \ref{p-4.3} to these chosen elements and passing to a subsequence gives integral-valued valuations $w_j$ centred at $o_j$, a constant $C>0$, and a positive integer $\ell'$, independent of $j$ in this subsequence.
Their restrictions $v_j$ to the meromorphic germ fields of $(X_{i_j},x_j)$ are integer valued.
As explained in Subsection \ref{s-5.2}, the $v_j$ are realized by divisors on compact modifications and hence are divisorial valuations over the fixed model $W$.
The discrepancy comparison in Subsection \ref{ss-5.1} and inequality \eqref{e-8} give
\[ \epsilon A_W(v_j) \leq A_0(v_j)\leq A_{i_j}(v_j) \leq A_{V_j}(w_j)\leq C,\]
where $\epsilon>0$ depends only on the original pair and resolution.
The character calculation in Subsection \ref{ss-5.1} also shows that
\(dA_{i_j}(v_j)\equiv_\Z \ell'a_j. \)
The constants $C$ and $\ell'$ may depend on $\alpha$, but they are fixed throughout the subsequence $i_j$.

Since $0\leq E_{i_j}(v_j)\leq A_{i_j}(v_j)\leq C$, we may pass to a further subsequence such that $E_{i_j}(v_j)\to e$.
The bound $A_W(v_j)\leq C/\epsilon$ allows us to apply Proposition \ref{lem:countable-limits}, which places $e$ in a countable set
$\mathcal E\subset\mathbb R$ fixed independently of $\alpha$.
Since $kdD_0$ is Cartier and the valuations $v_j$ are integer valued,
$kdA_0(v_j)\in\mathbb Z$ then \[
[kdE_{i_j}(v_j)]
=[kdA_{i_j}(v_j)-kdA_0(v_j)]
=[kdA_{i_j}(v_j)]
=[k\ell'a_j],\] where $[t]$ denotes the class of $t$ modulo $\mathbb Z$. Taking the limit in $\mathbb R/\mathbb Z$ therefore yields
\[  [kd e]=[k\ell'\alpha]. \]
For each fixed $e\in\mathcal E$ and $\ell'\in\mathbb Z_{>0}$, this equation allows only finitely many values of $\alpha\in[0,1)$.
Taking the union over the countable set $\mathcal E\times\mathbb Z_{>0}$ gives only countably many possible $\alpha$'s, contradicting the arbitrary choice of $\alpha$ in $[0,1)$.
Thus the local indices of all the $L_i$ have a common bound $r_0$. Then
\[  R:=2d\,\operatorname{lcm}(1,2,\ldots,r_0) \]
has the property that $RD_i$ is Cartier on every $X_i$, proving Proposition \ref{p-1}.

\vskip.3cm
\noindent{\bf AI statement.}
The main result of this paper was obtained with the aid of ChatGPT and is based on the  AI generated papers \cite{OAI1} and \cite{OAI2}. 
\section{Descending chain condition}
In this section we prove a result about the descending chain condition for "bounded" coherent ideals on a compact K\"ahler manifold.

Since $\mathcal O _X$ is noetherian, any ascending chain of coherent ideals $a_1\subset a_2\subset \ldots$ stabilizes, i.e., there exists an integer $j_0$ such that $a_j=\cup _{i>0}a_i$ for every $j\geq j_0$. This clearly fails for a descending sequence of ideals, e.g. $\mathfrak m\supset \mathfrak m^2\supset \mathfrak m^3\supset \ldots $ where $\frak m$ is a maximal ideal. However, under appropriate conditions, one can show that certain descending sequences stabilize. One example of this phenomenon is if $X$ is compact, $a_i$ is a descending sequence of coherent ideals and $L$ is a line bundle such that $L\otimes a_i$ is generated for all $i>0$. We then have $H^0(L\otimes a_1)\supset H^0(L\otimes a_i)\supset \ldots$ and since $H^0(L\otimes a_i)$ are finite dimensional vector spaces, then $H^0(L\otimes a_i)=H^0(L\otimes a_{i+1})$ for all $i\gg 0$. 
Since $L\otimes a_i$ is generated, then the inclusion $L\otimes a_i\supset L\otimes a_{i+1}$ is an equality for all $i\gg 0$ and hence $a_i=a_{i+1}$ for all $i\gg 0$.

In the analytic context we can not use cohomological methods and hence we rely on different positivity notions, such as nefness, instead of global generation.

\begin{lemma}\label{l-2.1} Let $W$ be a  compact
K\"ahler manifold of dimension $n \geq 1$. Let $\omega$ be a Kähler class and let $a_1 \supset 
a_2 \supset\ldots$ be nonzero integrally closed coherent ideals. Suppose that for every $i$
there is a log resolution $r_i : W_i \to W$ with $a_i\OO_{W_i} = \OO_{W_i} (-F_i)$ and $r_i^*\omega -[F_i]$ 
nef. Then the ideals stabilize.
\end{lemma}
\begin{proof} Let $r:W'\to W$ be a resolution of $(W, a_i)$, and  \[V(a_i):=\int _{W'}(2r^*\omega-[F'_i])^n,\qquad  a_i\OO_{W'}=\OO _{W'}(-F'_i).\]
It follows easily from the projection formula that $V(a_i)$ does not depend on the choice of $W'$, and since $r^*\omega -[F'_i]$ is nef, then $V(a_i)>0$. Suppose that $a_i\supset a_{i+1}$ is a proper inclusion and $r:W'\to W$ a common resolution of $a_i$ and $a_{i+1}$. Let $\rho :Z\to W$ be the normalized blow up of $a_{i+1}$ and $W'\to Z$ the induced morphism. Since $a_{i+1}$ is integrally closed, it is determined by its Rees valuations (see \cite[Theorem 4.1.6]{LJT08})  and since the inclusion $a_i\supset a_{i+1}$ is  proper, there is a divisor $\bar P$ on $Z$ such that ${\rm ord}_{\bar P}(a_i)<{\rm ord}_{\bar P}(a_{i+1})$.
Let $P$ be the strict transform of $\bar P $ on $W'$,
then $P\leq G=F'_{i+1}-F'_i$.  
Since $M_i=2r^*\omega -[F'_i]$ and $M_{i+1}=2r^*\omega -[F'_{i+1}]$ are both nef, we have
\[V(a_i)-V(a_{i+1})=\int _{W'}[G]\sum _{k=0}^{n-1}M_i^kM_{i+1}^{n-k-1}\geq \int _{P}M_{i+1}^{n-1}.\] We now let $S=r(P)$, $f$ a general fiber of $r|_P$, $d=\dim S$,  and $\bar f$ a general fiber of $\rho |_{\bar P}$. Note that $a_{i+1}\cdot \OO _Z=\OO _Z(-F_Z)$ where $-F_Z$ is ample over $W$. Then
\[\int _{P}M_{i+1}^{n-1}\geq \int _{P}M_{i+1}^{n-1-d}(r^*\omega ^d) =\int _{\bar P}(2\rho ^*\omega -[F_Z])^{n-1-d}(\rho ^*\omega ^d) =e\int _S\omega ^d \] where $e=\int _{\bar f} (2\rho ^*\omega -[F_Z])^{n-1-d}=\int _{\bar f} [-F_Z]^{n-1-d}>0$ as $-F_Z$ is relatively ample. Note that as $-F_Z$ is Cartier, then in fact $e\geq 1$.
Since $\omega $ is K\"ahler, there is a constant $C>0$ such that $\int _S\omega ^d\geq C$ for any $d$-dimensional analytic subvariety $S\subset W$, and any $1\leq d\leq n$. Replacing $C$ by ${\rm min}\{C,1\}$ we may assume that the stated inequality also holds when $\dim S=0$. Since $V(a_i)>0$, these decreases can happen at most $V(a_1)/C$ times and hence the ideals stabilize.
\end{proof}
\begin{lemma}\label{l-2.2} Let $(W, \omega )$ be a compact K\"ahler manifold, and $\mathcal U$ a family of quasi-plurisubharmonic functions with $dd^ cu \geq -C\omega$ for some constant $C \geq 0$. For each $s>0$, there
is a K\"ahler class $\Theta _s$ such that, for every $u \in  \mathcal U$ and every log resolution
$r : W'\to W$ of $I = \mathcal J (su)$ with $I\cdot \OO_{W'} = \OO_{W'} (-F )$, the class $r^*\Theta _s - [F ]$ is nef.
\end{lemma}
\begin{proof} 
We first prove that there is a constant $A_s>0$, depending only on
$(W,\omega)$, $C$, and $s$, such that, for every $u\in\mathcal U$, there is a
quasi-plurisubharmonic function $b_u$ on $W$ satisfying
\begin{equation}\label{eq:uniform-bergman-potential}  dd^c b_u\geq -A_s\omega, \qquad b_u=\log\left|\mathcal J(su)\right|+O(1).
\end{equation}

We choose finitely many coordinate balls
\[A_j\Subset B_j\Subset \Omega_j\]
such that the $W=\cup A_j$. After enlarging the coordinate balls if
necessary, we may choose smooth functions $q_j$ on $\Omega_j$ such that \[dd^c q_j\geq sC\omega.\]
Thus the function \[ \varphi_j:=su+q_j\]
is plurisubharmonic on $\Omega_j$.

Let $\mathcal H_j$ be the Hilbert space of holomorphic functions on
$\Omega_j$ with finite norm
\[\mathcal H_j:=\left\{  f\in \OO (\Omega _j)\ {\rm s.t.}\ \|f\|_j^2:=\int_{\Omega_j}|f|^2e^{-2\varphi_j}\,d\lambda_j<\infty\right\},\]
and define the associated Bergman kernel and Bergman potential by
\[ K_j(x)=\sup_{\|f\|_j\leq 1}|f(x)|^2, \qquad b_j=\frac12\log K_j.\]
Then $b_j$ is plurisubharmonic. Since \(\Omega_j\) is a bounded pseudoconvex coordinate ball, 
\cite[Proposition 5.7]{Dem12b} shows that \(\mathcal H_j\) generates \(\mathcal J(\varphi_j),\)
and we have
\begin{equation}\label{eq:local-singularity} b_j=\log|\mathcal J(\varphi_j)|+O(1) =\log|\mathcal J(su)|+O(1) \quad\text{on }B_j. 
\end{equation}

We claim that there is a constant $M_s>0$, independent of $u$, such that
\begin{equation}\label{eq:overlap}
 b_j\leq b_k+M_s,\qquad b_k\leq b_j+M_s
 \qquad\text{on }B_j\cap B_k \end{equation}
To see this, fix $x\in B_j\cap B_k$. By compactness, there is a radius
$\rho>0$, independent of $j,k,x$, such that the $j$-coordinate ball
$B_j(x,\rho)$ is contained in $\Omega_j\cap\Omega_k$. The
H\"ormander--Ohsawa--Takegoshi extension estimate, applied with the weight
\[ 2\varphi_j(z)+2n\log|z-x|+|z|^2,\]
gives a constant $L_s>0$, independent of $u,j,k,x$, such that for every
holomorphic function $f$ on $B_j(x,\rho)$ with \(\int_{B_j(x,\rho)}|f|^2e^{-2\varphi_j}\,d\lambda_j<\infty,
\) there exists $F\in \mathcal H _j$  satisfying
\[ F(x)=f(x),\qquad \int_{\Omega_j}|F|^2e^{-2\varphi_j}\,d\lambda_j\leq L_s\int_{B_j(x,\rho)}|f|^2e^{-2\varphi_j}\,d\lambda_j.\]
The functions $q_j-q_k$ and the Jacobians of the coordinate changes are
uniformly bounded on the finitely many overlaps. Consequently, restriction
of an element of $\mathcal H_k$ to $B_j(x,\rho)$, followed by the preceding
extension, yields
\[ K_j(x)\geq L_s^{-1}K_k(x).\]
Interchanging $j$ and $k$ gives the reverse inequality $L_sK_k(x)\geq K_j(x)$. This implies \eqref{eq:overlap}, after replacing $M_s$ by $\frac12\log L_s$.

Choose smooth functions $h_j$ on $B_j$ such that \[ 0\leq h_j\leq 1,\qquad  h_j=0 \ \text{on }A_j, \qquad h_j=1\ \text{near }\partial B_j.\]
Define
\[ b_u(x)=\max_{j:\,x\in B_j}\left\{b_j(x)-(M_s+1)h_j(x)\right\}.\]
This is globally well-defined. Indeed, if $x$ is close to the boundary of
$B_j$, then $x\in A_k$ for some $k$, then by \eqref{eq:overlap}
\[ b_j(x)-(M_s+1)h_j(x)\leq b_k(x)-1, \]
where \(h_j=1\) and \(h_k=0\). Thus, as $x$ approaches the boundary of $B_j$, $b_j(x)$ does not contribute to the maximum defining $b_u(x)$. Thus $b_u(x)$ is a globally well defined quasi psh function.

Since the functions $h_j$ are fixed and there are only finitely many of
them, there is a constant $A_s>0$, independent of $u$, such that \(dd^c b_u\geq -A_s\omega. \)
Furthermore, \eqref{eq:local-singularity} implies
\( b_u=\log\left|\mathcal J(su)\right|+O(1),\)
and this proves \eqref{eq:uniform-bergman-potential}.

Now set \(\Theta_s=(A_s+1)[\omega]\) and let
\(T_u=(A_s+1)\omega+dd^c b_u.\)
Then $T_u$ is a closed positive current representing $\Theta_s$.

Let $r:W'\to W$ be a log resolution of
$I:=\mathcal J(su)$, and write
\[ I\cdot\mathcal O_{ W'}=\mathcal O_{W'}(-F). \]
 If $g_F$ is a local equation of $F$, then
\(b_u\circ r=\log|g_F|+O(1).\)
Thus the generic Lelong number of \(r^*T_u\) along each component of
\(F\) equals its coefficient in \(F\). By  Siu decomposition,
\(r^*T_u=[F]+R_u\) for a closed positive current \(R_u\). Since
\(b_u\circ r-\log|g_F|\) is locally bounded, then 
\(R_u=r^*T_u-[F]\) has locally bounded potentials. In particular, \(R_u\)
has vanishing Lelong numbers.
Its cohomology class
is
\[ \bigl\{r^*T_u-[F]\bigr\}=r^*\Theta_s-[F].\]
By Demailly's nefness criterion for a class represented by a positive
current with zero Lelong numbers, the class $r^*\Theta_s-[F]$ is nef (see \cite[Proposition 17.8(a)]{Dem12b}).
\end{proof}
\begin{proposition}\label{p-2.3}Let $(W, \omega )$ be a compact K\"ahler manifold, and $\mathcal U$ a family of quasi-plurisubharmonic functions with $dd^ cu \geq -C\omega$ for some constant $C \geq 0$. Suppose that $u_1\geq u_2,\ldots $ is a decreasing sequence in $\mathcal U$, then, for any  $s>0$ the ideals $\mathcal J(su_i)$ stabilize and so $\mathcal J(su_i)=\cap _{j>0}\mathcal J(su_j)$ for any $i\gg 0$.  \end{proposition}
\begin{proof} We have a decreasing sequence $\mathcal J(su_1)\supset \mathcal J(su_2)\supset \ldots$ of integrally closed ideals. The claim follows from Lemmas \ref{l-2.1} and \ref{l-2.2}.
\end{proof}

\section{Discrepancy bounds}
Recall the notation from Subsection \ref{ss-sketch}.  In particular we assume that $K_X+B$ is pseudo-effective, $\omega$ is K\"ahler, $K_{X_i}+B_i+t\omega _i$ is nef for $\lambda _i\geq t\geq \lambda_{i+1}$, $\lambda _i>0$, and that $\lim _{i\to \infty}\lambda _i=0$. Since there are at most finitely many divisorial contractions, we may also assume that for $i\gg 0$, each $X_i\dasharrow X_{i+1}$ is a flip.
If $v$ is the valuation corresponding to a divisor over $X$, then $A_i(v)$ denotes the log discrepancy with
respect to $(X_i,B_i)$. We define \[E_i(v)=A_i(v)-A_0(v),\qquad \Delta _i(v)=A_{i+1}(v)-A_i(v).\]
Thus $E_i(v)=\sum _{j=0}^{i-1}\Delta _j (v)$. Note that $\Delta _i(v)\geq 0$ as $\psi_i:X_i\dasharrow X_{i+1}$ is a $K_{X_i}+B_i$ flip or divisorial contraction.  If $\psi _i$ is a flip, then $f_i:X_i\to Z_i$ is the flipping contraction and $f_i^+:X_{i+1}\to Z_i$ is the flipped contraction. If $\psi _i$ is a divisorial contraction, then $f_i=\psi _i$, $Z_i=X_{i+1}$, and $f_{i}^+={\rm id}_{Z_i}$.
Since $Z_i$ has rational singularities and $K_{X_i}+B_i+\lambda _{i+1}\omega _i\equiv _{Z_i}0$, then there is a class $\zeta _i\in H^{1,1}_{\rm BC}(Z_i)$ such that
\[K_{X_i}+B_i+\lambda _{i+1}\omega _i=f_i^*\zeta _i,\qquad K_{X_{i+1}}+B_{i+1}+\lambda _{i+1}\omega _{i+1}=(f_i^+)^*\zeta _i,\]
see \cite[Lemma 2.39 and Proposition 4.2]{HX26}.
On a common projective resolution $q_i:Y\to X_i$ and $q_{i+1}:Y\to X_{i+1}$ we let
\[\Delta _i^Y=q_i^*(K_{X_i}+B_i)-q_{i+1}^*(K_{X_{i+1}}+B_{i+1})\geq 0.\]
For any divisorial valuation $v$ over $Y$ we have $v(\Delta _i^Y)=\Delta _i(v)$. In other words, $\Delta _i$ is the b-divisor defined by the valuations $\Delta _i(v)$ and $\Delta _i^Y$ is the trace on $Y$.
Let $D=K_{X}+B$  and $D_i=K_{X_i}+B_i$. We have \[    D_{i}+t\omega _i\equiv _{Z_i}(1-\frac t{\lambda _{i+1}})(K_{X_i}+B_i)\qquad  D_{i+1}+t\omega _{i+1}\equiv _{Z_{i}}(1-\frac t{\lambda _{i+1}})(K_{X_{i+1}}+B_{i+1}),\] 
\[q_i^*([D_{i}]+t\omega _i)-q_{i+1}^*([D_{i+1}]+t\omega _{i+1})= (1-\frac t{\lambda _{i+1}})[\Delta ^Y_i]\qquad {\rm in\ } H^{1,1}_{\mathrm{BC}}(Y,\mathbb R).\]

Fix $p:W\to X$ a projective log resolution of $(X,B)$. Choose a smooth closed real (1, 1)-form $\theta$ representing $[p^*(K_X+B)]$, and set $\gamma :=p^*\omega$.
In particular $\gamma$ is smooth and semipositive. For $t>0$, define \(u_t\) and \(T_t\) by
\[u_t:=P_{\theta+t\gamma}(0)=
\left( \sup\{u\in\operatorname{PSH}(W,\theta+t\gamma):u\leq0\}\right)^*, \qquad T_t:=\theta+t\gamma+dd^cu_t.\]  Here the superscript \(^*\) means upper semicontinuous regularization. Then \(u_t\) is the largest \((\theta+t\gamma)\)-plurisubharmonic function on \(W\) that is everywhere $\leq 0$ and \(T_t\) is a positive current with minimal singularities in \([p^*D]+t[\gamma]\).
Since \(\gamma\geq0\), the potentials \(u_t\) decrease as \(t\to 0^+\).
Let \( \sigma_v(t):=\nu_v(u_t).\)
 Choose
$\psi\in\operatorname{PSH}(W,\theta)$ with $\sup_W\psi=0$.
Then
\(\psi\leq u_t\leq0.\)
Semipositivity of $\gamma$ implies $u_{t'}\leq u_t$ whenever $0<t'\leq t$.
For a fixed K\"ahler form $\omega_W$ on $W$, we also have
\[ dd^cu_t\geq-\theta-t\gamma\geq-C\omega_W
 \qquad(0<t\leq1) \]
with $C$ independent of $t$.
For every divisorial valuation $v$,
\[  0\leq\sigma_v(t)\leq\nu_v(\psi)<\infty. \]
Consequently, the limit
\(  \sigma_v:=\lim_{t\to 0^+}\sigma_v(t) \)
exists.

\begin{proposition}\label{p-5} For every
 $s>0$ there exists a nonzero coherent ideal $\mathcal I_s$ on $W$ and an integer
$i(s)$ such that for every $i \geq  i(s)$ and every divisorial valuation $v$, we have
\[\frac{v(\mathcal I_s)}s\leq E_i(v)\leq \frac{v(\mathcal I_s)+A_W(v)}s\] where $A_W(\ldots )$ denotes the log discrepancy with respect to $(W,0)$. In particular for every $v$,
\[0\leq \Delta _i(v)=E_{i+1}(v)-E_i(v)\leq \frac 1 s A_W(v) \qquad {\rm for\ all\ } i\gg 0.\]
\end{proposition}
\begin{proof} Let $q_0:Y\to X$ be a common projective log resolution of $(X_j,B_j)$ for $0\leq j\leq i$ and $q_j:Y\to X_j$ the induced morphisms. 
We have
\[q_0^*([D]+t\omega) =q_i^*([D_i]+t\omega_i)+[F_{i,t}], \qquad
 F_{i,t}:=\sum_{j=0}^{i-1}      \left(1-\frac{t}{\lambda_{j+1}}\right)\Delta_j^Y.\]
The divisor $F_{i,t}$ is effective and
$q_i$-exceptional, for $0<t\leq\lambda_i$. 
Since $[D_i]+\lambda_i\omega_i$ is nef, then by \cite[Lemma A.5]{DHY26}, we have 
\[ N\bigl(q_0^*([D]+\lambda_i\omega)\bigr)
 =N\bigl(q_i^*([D_i]+\lambda_i\omega_i)\bigr)+F_{i,\lambda_i}= F_{i,\lambda_i}.\]
For any given divisorial valuation $v$, we may choose $Y$ extracting $v$.
It follows that
\[ \sigma_v(\lambda_i)= \sum_{j=0}^{i-1} \left(1-\frac{\lambda_i}{\lambda_{j+1}}\right)\Delta_j(v) \leq \sum_{j=0}^{i-1} \Delta_j(v) =E_i(v).\]
For fixed $i$ and $k\geq i$, we have
\[ \sigma_v(\lambda_k) \geq\sum_{j=0}^{i-1} \left(1-\frac{\lambda_k}{\lambda_{j+1}}\right)\Delta_j(v). \]
Letting $k\to\infty$ yields $E_i(v)\leq\sigma_v$. Thus
\begin{equation}\label{e-1}   \sigma_v(\lambda_i)  =  \sum_{j=0}^{i-1}  \left(1-\frac{\lambda_i}{\lambda_{j+1}}\right)\Delta_j(v)  \leq E_i(v)\leq\sigma_v. \end{equation}
Fix $s>0$. The functions $su_{\lambda_i}$ decrease as $i\to \infty$, and have a uniform
lower curvature bound. Proposition \ref{p-2.3}, the sequence $\mathcal J(su_{\lambda _i})\supset \mathcal J(su_{\lambda _{i+1}})$ stabilizes for $i\gg 0$. Therefore, we have a coherent ideal
$\mathcal I_s$ and an integer $i(s)$ such that
\(\mathcal J(su_{\lambda_i})=\mathcal I_s\)
for \(i\geq i(s)\).
By \cite[Theorem 5.1 and (5.3)]{BFJ08}, it follows that
\begin{equation}\label{e-2} s\sigma_v(t)-A_W(v) \leq v\bigl(\mathcal J(su_t)\bigr) \leq s\sigma_v(t). \end{equation}
For $i\geq i(s)$ and $t=\lambda _i$, the right-hand inequality gives
$v(\mathcal I_s)/s\leq\sigma_v(\lambda_i)$.
Letting $i\to\infty$ (and hence $t\to 0^+$) in the left-hand inequality gives
\begin{equation}\label{e-3}  \sigma_v\leq\frac{v(\mathcal I_s)+A_W(v)}s. \end{equation}
Combining \eqref{e-3} and \eqref{e-1}, one sees that
\begin{equation}\label{e-4} \frac{v(\mathcal I_s)}s  \leq E_i(v)\leq\sigma_v  \leq\frac{v(\mathcal I_s)+A_W(v)}s. \end{equation}
Applying these bounds to $E_i(v)$ and $E_{i+1}(v)$ proves the final
assertion. The index $i(s)$ is independent of $v$.

\end{proof}
We write $K_W+B_W=p^*(K_X+B)$ and $S_W={\rm Supp}(B_W)$. Since $(X,B)$ is klt, every coefficient of \(B_W\) is \(<1\), and so there exists a constant $0<\epsilon <1$ such that $(1-\epsilon)S_W\geq B_W$. Since $(W,S_W)$ is log smooth and log canonical, then $v(S_W)\leq A_W(v)$ and so
\begin{equation}\label{e-8}A_i(v)\geq A_0(v)=A_W(v)-v(B_W)\geq \epsilon A_W(v).
 \end{equation}
\begin{proposition}\label{lem:countable-limits}
There is a countable subset $\mathcal E\subset \mathbb R$ with the following
property. For any subsequence $i_j\to\infty$, and  integer-valued divisorial
valuations $v_j$ over $W$ such that
\[ {\sup_{j>0}} A_W(v_j)<\infty , \]
 if the sequence $E_{i_j}(v_j)$ converges to $e=\lim E_{i_j}(v_j)$, then $e\in \mathcal E$.
\end{proposition}

\begin{proof} We let $\Sigma$ be the divisor on $W$ whose components are either $p:W\to X$ exceptional or the strict transforms of the finitely many divisors contracted by the given mmp.
If a prime divisor $Q\subset W$ is not a component of $\Sigma$, then it is not contracted by $W\dasharrow X_i$, and therefore
\(E_i(\operatorname{ord}_Q)=0.\)
By Proposition \ref{p-5} we have $\operatorname{ord}_Q(\mathcal I_s)=0$ for every
integer $s\geq1$. In particular,
\({\rm Supp}({\rm div}_{W}(\mathcal I_s))
 \subseteq\ \Sigma\).
We now pick a tower of projective modifications
\[ \cdots\longrightarrow W_s\longrightarrow W_{s-1} \longrightarrow\cdots\longrightarrow W_0=W\]
such that each $W_s$ is smooth and $\pi_s:W_s\to W$ is a common log resolution of
$(W,\Sigma)$, $(W,\mathcal I_1),\ldots,(W,\mathcal I_s)$.
Note that the $W_j$ are chosen independently from the sequences $i_j,v_j$.

Fix an integer $N\geq1$ and a sequence as in the statement with
$A_W(v_j)\leq N$.

\smallskip
\noindent{\bf Step 1.}  
For each fixed $s\geq1$, we let
\[ b_{s,j}:=v_j(K_{W_s/W})=A_W(v_j)-A_{W_s}(v_j).\]
Note that these are bounded integers $0\leq b_{s,j}\leq N$. Therefore, passing to a diagonal subsequence, we may assume that $b_{s,j}$ is eventually constant in $j$
for every fixed $s$. We let \[b_s=\lim b_{s,j}=b_{s,k}\qquad {\rm for}\ k\gg 0.\]
Let $h_s:W_{s+1}\to W_s$, then 
\[ K_{W_{s+1}/W}=K_{W_{s+1}/W_s}+h_s^*K_{W_s/W},\] where $K_{W_{s+1}/W_s}\geq {\rm div.Ex.}(h_s)$.
This implies that $b_{s+1}\geq b_s$, and hence there is an $M\geq1$ such that
$b_t=b_M$ for all $t\geq M$.
Consequently, for every fixed $t\geq M$,
\[  \lim_{j\to\infty}v_j(K_{W_t/W_M})= \lim_{j\to\infty}(b_{t,j}-b_{M,j})=0.\]

\smallskip
\noindent{\bf Step 2.} Let $\mathcal Q_M$ be the finite set consisting of the exceptional prime
divisors of $\pi_M$ and the strict transforms on $W_M$ of the components
of $\Sigma $.
For each $Q\in\mathcal Q_M$, the pair $(W_M,Q)$ is log canonical, so
\[0\leq v_j(Q)\leq A_{W_M}(v_j)\leq A_W(v_j)\leq N.\]
After passing to a further subsequence, each of these integer values is eventually
constant so that
\[ v_j(Q)=c_Q,\qquad Q\in\mathcal Q_M,\ c_Q\in\{0,\ldots,N\}. \]

\smallskip
\noindent{\bf Step 3.} For every integer $s\geq1$, we let
\[G_{M,s}:= {\rm div}_{W_M}(\mathcal I_s\cdot \OO_{W_M}) = \sum_{Q\in\mathcal Q_M}{\rm ord}_Q(\mathcal I_s)\,Q\]
be the divisorial part of the ideal on $W_M$.
Choose $t\geq\max\{M,s\}$ and let $h:W_t\to W_M$.
Since $t\geq s$, then $W_t\to W$ is a log resolution of $\mathcal I_s$, and so we may write
\[ \mathcal I_s\cdot \OO_{W_t}=\OO_{W_t}(-F_{s,t})\] for some effective divisor $F_{s,t}$.
The residual ideal
\[\mathfrak a_{M,s}:=  (\mathcal I_s\cdot \mathcal O_{W_M})\otimes\mathcal O_{W_M}(G_{M,s})\subseteq\mathcal O_{W_M}\]
has cosupport of codimension $\geq 2$.
Its pullback to $W_t$ is invertible, and we may write
\[F_{s,t}=h^*G_{M,s}+R_{s,t}, \]
where $R_{s,t}$ is effective and $h$-exceptional.
Every $h$-exceptional prime divisor has strictly positive coefficient in
$K_{W_t/W_M}$, so there is a constant $C_{s,t}>0$ with
\[ 0\leq R_{s,t}\leq C_{s,t}K_{W_t/W_M}.\]
By Step~1 we have $\lim _{j\to \infty}v_j(K_{W_t/W_M})= 0$ and hence $\lim _{j\to \infty}v_j(R_{s,t})= 0$. Thus the following limit exists 
\[a_s:=\lim_{j\to\infty}v_j(\mathcal I_s) = \sum_{Q\in\mathcal Q_M}\operatorname{ord}_Q(\mathcal I_s)c_Q.\]
This formula holds for every $s\geq1$, including $s<M$.
Hence, the full sequence $(a_s)_{s\geq1}$, is determined by $M$ and the
finite set of integers $(c_Q)_{Q\in\mathcal Q_M}$.

\smallskip
\noindent{\bf Step 4.} Suppose $E_{i_j}(v_j)\to e$.
For every fixed positive integer $s$, Proposition \ref{p-5} applies once
$i_j\geq i(s)$, and passing to the limit as $j\to \infty$ gives
\[\frac{a_s}s\leq e\leq\frac{a_s+N}s.\]
For arbitrary integers $N,M\geq1$ and
$c\in\{0,\ldots,N\}^{\mathcal Q_M}$, define
\[a_s(M,c):= \sum_{Q\in\mathcal Q_M}\operatorname{ord}_Q(\mathcal I_s)c_Q, \qquad \mathcal E_{N,M,c}:= \bigcap_{s\geq1} \left[\frac{a_s(M,c)}s,\frac{a_s(M,c)+N}s\right].\]
Each $\mathcal E_{N,M,c}$ is either empty or a singleton, since the interval
lengths $N/s$ can be chosen arbitrarily small. The preceding argument shows that every limit in
the statement belongs to
\[ \mathcal E:=
 \bigcup _{N, M\geq1}\left( \bigcup _{c\in\{0,\ldots,N\}^{\mathcal Q_M}}\mathcal E_{N,M,c} \right).\]
For each $N,M$ there are only finitely many such tuples $c\in\{0,\ldots,N\}^{\mathcal Q_M}$, and therefore
$\mathcal E$ is countable.
\end{proof}

\section{Orbifold and analytic quotients by finite groups}
We will need the following result see \cite[Theorem 3.1]{OAI1}.
\begin{theorem}\label{t-4.1} Let $(V_i, o_i)$ be a sequence of pointed algebraic Gorenstein
klt germs of a fixed dimension $n \geq 2$. Suppose that a finite cyclic group $G_i$ acts on
an affine representative of each germ, fixes $o_i$, and admits a semi-invariant canonical generator $\theta _i$,
a volume form with zero divisor near $o_i$. Choose an arbitrary element $\gamma _i\in G_i$ for each $i$.
After passage to a subsequence, there are divisorial valuations $w'_i$ centred at $o_i$, a constant $C > 0$,
and a positive integer $l$, independent of $i $ on that subsequence, such that
\[A_{\theta _i} (w'_
i) \leq C\]
and, for every nonzero rational semi-invariant function $h$ on $V_i$ of character $\chi$,
 \[w'_i(h) \equiv l \ {\rm angle}(\chi(\gamma _i)) \qquad ({\rm mod}\ \Z).\]
In particular, $w'_i$ restricts to an integral-valued valuation of $\C(V_i)^{G_i}$.
\end{theorem}
We record the following convention. For \(\zeta\in S^1\), we let \(\operatorname{angle}(\zeta)\in\mathbb R/\mathbb Z\) be the class determined by \(\zeta=e^{2\pi {\rm i}a}\), and \(\arg\zeta=a\) denotes its representative in \([0,1)\).
The following result allows us to pass from the analytic setting to the algebraic setting.
\begin{proposition}\label{p-4.2} 
Let $(V,o)$ be a Gorenstein klt complex analytic germ and $G$ a
finite group acting on $V$ and fixing $o$.  Then there exist a $G$-equivariant projective log
resolution $f\colon Y\to V$, an algebraic Gorenstein klt germ $(V',o')$
with a $G$-action fixing $o'$, and a $G$-equivariant projective log resolution
$f'\colon Y'\to V'$ with the following properties.  The scheme-theoretic fibres $f^{-1}(o)$ and $(f')^{-1}(o')$ are
$G$-equivariantly identified, preserving the marked exceptional divisors.  At
corresponding points $y$ and $y'$ of these fibres we have the following.

\begin{enumerate}
\item The representations of the stabilizer groups on \(T_yY\) and \(T_{y'}Y'\) are isomorphic.
\item The correspondence between exceptional divisors through $y$ and $y'$
is $G_y$-equivariant. For corresponding exceptional divisors $E_j,E'_j$, put
\[
 G_{y,E_j}:=\{g\in G:g(y)=y,\ g(E_j)=E_j\}.
\]
The markings identify this subgroup with $G_{y',E'_j}$, and the normal
lines $N_{E_j/Y,y}$ and $N_{E'_j/Y',y'}$ give isomorphic representations
of this subgroup.
\item The discrepancy coefficients of corresponding exceptional divisors in $K_Y-f^*K_V$ and $K_{Y'}-(f')^*K_{V'}$ agree.
\end{enumerate}

If the action of $G$ on $(V,o)$ is faithful, then the action of $G$ on $(V',o')$ may be chosen to be
faithful.
\end{proposition}
\begin{proof}
If \(V\) is smooth, then the action of the finite group \(G\) is locally
analytically linearizable. Thus we may take \(V'=T_oV\), with its linear
\(G\)-action, and take both resolutions to be the identity. Hence assume
that \(V\) is singular, and write \(n=\dim V\).

\medskip
\noindent{\bf Step 1.}
Choose a \(G\)-equivariant embedding of \((V,o)\) into a smooth representation germ, and write
\[  S=\mathbb C\{z_1,\ldots,z_N\},\qquad S_0=\mathbb C[z_1,\ldots,z_N]_{(z)},\qquad
 R=S/I, \]
where \((V,o)=\operatorname{Specan}(R)\), and \(G\) acts linearly on the coordinates, compatibly on
\(S_0\) and \(S\) and $I$ is $G$ stable.

We first apply functorial resolution and principalization to the
excellent scheme ${\rm Spec}R$, and then analytify near its
closed fibre (see \cite[Theorem 1.1.11]{Tem}). Functoriality ensures that the construction is $G$-equivariant.
The projective resolution, its marked Cartier divisors, and their
line bundles and sections are defined by finite-presentation data over $R$.
 We denote the \(G\)-equivariant projective log
resolution by \(f:Y\to V\). It is an isomorphism over \(V_{\rm reg}\), with
reduced exceptional divisor
\(E=\sum_j E_j.\)
After replacing the resulting relatively ample exceptional divisor by a sufficiently
positive integral combination of the exceptional divisors of the successive
blow-ups, we may write
\[ Y=\operatorname{Bl}_J V,\qquad J\cdot \mathcal O_Y=\mathcal O_Y(-\Delta),\qquad \Delta=\sum_j m_jE_j,\quad m_j>0,\]
where \(J\subset R\) is \(G\)-invariant, its cosupport has codimension at
least two, and \(-\Delta\) is \(f\)-ample and $G$-invariant. After taking a sufficiently
high Veronese we may assume that \(Y={\rm Proj}_V \oplus_{k\geq 0} f_*\OO _Y (-k\Delta)\).

\medskip
\noindent{\bf Step 2.}
Let
\(T_0=S_0^G\) and \(T=S^G \), then \(\widehat T_0\simeq\widehat T\). The natural
map
\[ S_0\otimes_{T_0}T\longrightarrow S \]
 is a map of finite \(T\)-modules that becomes an
isomorphism after completion and therefore it is an isomorphism. Since \(T_0\) is excellent,
\(T_0\to\widehat T_0\simeq\widehat T\) is regular. As
\(T\to\widehat T\) is faithfully flat, regularity descends to
\(T_0\to T\). By Popescu's theorem (\cite[Tag 07GC]{Stacks}),
\[T=\varinjlim_\lambda B_\lambda,\qquad \lambda \in \Lambda\]
with \(B_\lambda\) smooth over \(T_0\).
Pick a sufficiently large \(\lambda \in \Lambda\), and put
\[B=B_\lambda,\qquad P=S_0\otimes_{T_0}B.\]
The group \(G\) acts on \(P\) through the first factor. If \(\lambda \in \Lambda\) is sufficiently large, then we may  descend the
following finite data.
\begin{enumerate}
\item A $G$-stable generating set of $I$ and the identities expressing its
$G$-invariance, together with the matrices of a finite free $S$-resolution
of $R$ whose degree-zero term is $S$ and whose first differential has
these ideal generators as entries.
\item A finite $G$-equivariantly indexed generating set of $J$, a finite
homogeneous presentation of its Rees algebra, and the identities obtained
by substituting the ideal generators into the homogeneous equations.
\item The indexed family of line bundles and sections defining the
marked divisors $E_j$, the isomorphisms by which $G$ permutes this family,
their cocycle identities, and the factorization
\[  \mathcal O_Y^{\oplus e}\longrightarrow \mathcal O_Y(-\Delta) \xrightarrow{\ \prod_js_j^{m_j}\ }\mathcal O_Y\]
of the generator map for $J\cdot \mathcal O_Y$.
\end{enumerate}
In (1) we include the finitely many equations saying that successive
differentials compose to zero. Thus the descended matrices form a
complex. We will prove its exactness using the analytic graph below.
All the indicated objects, morphisms, and equations are of finite
presentation and descend for sufficiently large \(\lambda\in\Lambda\).

Let \(I_P\subset P\) be the descended ideal, and set
\[ A_P=P/I_P,\qquad V_P=\operatorname{Spec}A_P.\]
Let \(J_P\subset A_P\) be the descended ideal, and let
\(f_P:W_P\longrightarrow V_P\)
be the projective morphism defined by the descended homogeneous Rees
equations. After base change along $P \to S$, these schemes and the marked data recover
the original ones over $R$. At this point it is not clear that \(W_P\) is the blow-up of
\(J_P\). The evaluation identities in (2) give a closed immersion
\[
 \operatorname{Bl}_{J_P}V_P\hookrightarrow W_P.
\]

Let $p\in V_P$ be the image of $o$, and let $\widetilde p$ be its image
in $U_P:=\operatorname{Spec}P$. These points are fixed by $G$.
After localizing $B$ near the image of $p$, choose relative
\'etale coordinates $t_1,\ldots,t_b$, centred at that point.
The analytic ambient local ring is
\[Q:=\mathcal O_{U_P^{\mathrm{an}},\widetilde p}\simeq\mathbb C\{z_1,\ldots,z_N,t_1,\ldots,t_b\}. \]
The map $P\to S$ sends $t_j$ to a $G$-invariant convergent function $a_j(z)$, and the graph equations
\[  u_j=t_j-a_j(z),\qquad 1\leq j\leq b, \]
are $G$-invariant. They cut out the original germ $V$ in
$V_P^{\mathrm{an}}$.

Put $\Lambda=\mathbb C\{u_1,\ldots,u_b\}$. The coordinate change gives
$Q\simeq\mathbb C\{z,u\}$, so $Q$ is flat over $\Lambda$.
The descended free complex, extended to $Q$, consists of finite
$Q$-modules which are $\Lambda$-flat. Its reduction modulo $(u)$ is
the chosen free $S$-resolution of $R$. By \cite[Tag~00MI]{Stacks}, this
complex is exact in positive degrees and its degree-zero cokernel
\[  A:=Q/I_PQ=\mathcal O_{V_P^{\mathrm{an}},p}\]
is flat over $\Lambda$. Therefore $u_1,\ldots,u_b$ is an $A$-regular
sequence and
\[ A/(u_1,\ldots,u_b)\simeq R.\]
Since $R$ is Gorenstein and $u_1,\ldots,u_b$ is a regular sequence, then $A$ is Gorenstein
of dimension $n+b$.

We next check normality. After shrinking near $p$, the graph equations
remain a regular sequence along $V$. At every point of
$V_{\mathrm{reg}}$, regularity of the quotient by that sequence implies
regularity of $V_P^{\mathrm{an}}$. If its singular locus had a component
through $p$ of codimension at most one, intersection with the graph,
which is cut out by $b$ equations, would have dimension at least $n-1$.
This intersection can not be contained in $\operatorname{Sing}(V)$,
whose dimension is at most $n-2$. This contradicts regularity of the
total space at the regular points of the graph. Hence $A$ is
$R_1$, and since it is Gorenstein, it is $S_2$ and hence normal by Serre's criterion.

The algebraic and analytic local rings at $p$ have the same completion.
Using excellence and the completion criteria for normality and
Gorensteinness, the algebraic germ $(V_P,p)$ is normal and Gorenstein
as well. We may therefore shrink equivariantly around $p$ so that
$V_P$ is integral, normal, and Gorenstein of dimension $n+b$.
In what follows, analytic graph equations in algebraic completed local
rings mean their images under the identification with the corresponding
analytic completions.

\medskip
\noindent{\bf Step 3.}
Write $\widetilde W=\operatorname{Bl}_{J_P}V_P$.
After pulling back to the analytic graph $V$ we have closed immersions
\[ Y=\operatorname{Bl}_JV \hookrightarrow \widetilde W^{\mathrm{an}}\times_{V_P^{\mathrm{an}}}V \hookrightarrow W_P^{\mathrm{an}}\times_{V_P^{\mathrm{an}}}V=Y.\]
Their composition is the identity on the homogeneous Rees
coordinates, so both immersions are isomorphisms.
In particular, the central fibre
\[F:=f^{-1}(o)=(f_P)^{-1}(p)\]
is identified scheme-theoretically, and every closed point $y\in F$
has a corresponding point $\widetilde y\in\widetilde W$.

Set $C_y=\mathcal O_{W_P,y}$ and
$\widetilde C_y=\mathcal O_{\widetilde W,\widetilde y}$.
The nonzero ideal $J_P$ lies in the integral germ $V_P$, so
$\widetilde W$ is integral and birational over $V_P$. Consequently
\(\dim\widetilde C_y=n+b.\)
The closed immersion gives a surjection $C_y\twoheadrightarrow
\widetilde C_y$. On the other hand,
\[\widehat C_y/(u_1,\ldots,u_b) \simeq\widehat{\mathcal O}_{Y,y},\]
and $\widehat{\mathcal O}_{Y,y}$ is regular of dimension $n$.
Lifting its $n$ parameters and adjoining the $u_j$ generates the maximal
ideal of $\widehat C_y$. Thus
\[n+b\leq\dim\widehat C_y\leq\operatorname{edim}\widehat C_y\leq n+b. \]
It follows that $\widehat C_y$ is regular. The kernel of the surjection 
$\widehat C_y\twoheadrightarrow\widehat{\widetilde C_y}$ is zero, since
a nonzero ideal in the regular local domain $\widehat C_y$ would lower
the dimension.
Faithful flatness of completions shows that
$\widetilde W\hookrightarrow W_P$ is an isomorphism near every point
of $F$, and $W_P$ is smooth there. Since $f_P$ is proper, we may further shrink $V_P$ around $p$ 
($G$-invariantly) so that
\[W_P=\operatorname{Bl}_{J_P}V_P\quad\text{and}\quad W_P\text{ is smooth}.\]

Let $s_j$ now denote the descended marked sections. The descended map
from the generators of $J_P$ to the invertible sheaf corresponding to
$\mathcal O_Y(-\Delta)$ restricts to a surjection on $Y$.
Its cokernel modulo the graph ideal is therefore zero at every point
of $F$. Nakayama's lemma gives surjectivity near $F$.

The graph equations have independent differentials along $F$, since
$W_P$ and $Y$ are smooth of dimensions $n+b$ and $n$.
At a point of $F$, the restrictions of the marked sections vanishing
there define the SNC divisor $E$ on $Y$. Their differentials on $Y$
are independent. Together with $du_1,\ldots,du_b$, their differentials
are therefore independent on $W_P$. After a further shrinking, the
zero divisors $D_j=(s_j=0)$ are effective Cartier divisors with simple
normal crossings, and $Y$ is transverse to every marked stratum.
The descended factorization and its surjective first map give
\[  J_P\cdot\mathcal O_{W_P} =\mathcal O_{W_P}\!\left(-\sum_jm_jD_j\right). \]
All the fibre identifications, markings, and maps in this construction
are $G$-equivariant.

\medskip
\noindent{\bf Step 4.}
Choose an integer $q\geq2$,  and let $a'_j(z)$
be the Taylor polynomial of $a_j(z)$ of degree less than $q$.
Since the action of $G$ on the $z$-coordinates is linear and $a_j$ is
$G$-invariant, $a'_j$ is also $G$-invariant. Set
\[ u'_j=t_j-a'_j(z),\qquad V'=\mathrm V(u'_1,\ldots,u'_b)\subset V_P,\qquad Y'=W_P\times_{V_P}V',\]
and let $f':Y'\to V'$ be the induced morphism and $o'=p$.
These are algebraic objects with their induced $G$-actions.
Both graphs contain $p$, and their scheme-theoretic central fibres are
canonically identified:
\[ F=Y\times_V\{o\} =W_P\times_{V_P}\{p\}
   =Y'\times_{V'}\{o'\}.\]
These identifications are $G$-equivariant.

Working in the analytic local rings at points of the common fibre \(F\), we first compare \(Y\) and \(Y'\), together with the restrictions of the marked divisors \(D_j\).
For every $y\in F$,
\[  f_P^*(u'_j-u_j)  \in  \bigl(\mathfrak m_p\mathcal O_{W_P^{\mathrm{an}},y}\bigr)^q
 \subseteq \mathfrak m_y^2. \]
Thus the differentials of the two ordered systems of graph equations
agree at $y$.
The original graph has independent differentials and is transverse to
every stratum of the marked simple normal crossing divisor.
The same is therefore true of the new graph.
In particular, $Y'$ is smooth of dimension $n$ along $F$, and the
restrictions $E'_j=D_j|_{Y'}$ form a simple normal crossing divisor
there.
Under the identification of corresponding fibre points, we have
\[T_yY=T_yY'\subset T_yW_P\]
as $G_y$-representations. For each marked divisor \(D_j\) through \(y\), the divisors \(E_j=D_j|_Y\) and \(E'_j=D_j|_{Y'}\) have tangent spaces
\[T_yE_j=T_yY\cap T_yD_j=T_yY'\cap T_yD_j=T_yE'_j.\]
Here \(Y'\) and \(D_j\) are viewed analytically, and the corresponding fibre points are identified with the same point \(y\in F\).
Consequently the normal lines to $E_j$ and $E'_j$ agree as representations of
\[ G_{y,E_j}=\{g\in G:g(y)=y,\ g(E_j)=E_j\}.\]
The correspondence between the marked divisors and their normal lines is \(G_y\)-equivariant, including the permutation of the marked divisors through \(y\).
Moreover, every scheme-theoretic intersection of marked divisors with
$F$ is obtained by restricting the same sections defining the $D_j$.
Thus all these marked intersections are preserved as schemes with
their group actions.

We will next choose $q$ large enough to preserve bounds on the dimension of $V$ and on the the dimension of the cosupport of $J_P$.
We start by making the following elementary observation.
Let $(B,\mathfrak m)$ be a Noetherian local ring, let
$v_1,\ldots,v_b\in\mathfrak m$, and suppose that
$\dim B/(v_1,\ldots,v_b)\leq k$.
Choose $x_1,\ldots,x_k\in\mathfrak m$ and $M\geq1$ such that
\[  \mathfrak m^M  \subseteq(v_1 \ldots,v_b,x_1,\ldots,x_k). \]
If $v'_j-v_j\in\mathfrak m^{M+1}$ and
$I'=(v'_1,\ldots,v'_b,x_1,\ldots,x_k)$, then
\[ \mathfrak m^M\subseteq I'+\mathfrak m^{M+1}. \]
The finite $B$-module $Q_0=(\mathfrak m^M+I')/I'$ therefore satisfies
$Q_0=\mathfrak mQ_0$. Nakayama's lemma gives $Q_0=0$, and hence
$\mathfrak m^M\subseteq I'$.
Passing to \(\overline B=B/(v'_1,\ldots,v'_b)\), we obtain
\[\overline{\mathfrak m}^{\,M}\subseteq(\overline x_1,\ldots,\overline x_k).\]
Hence \((\overline x_1,\ldots,\overline x_k)\) is \(\overline{\mathfrak m}\)-primary, and Krull’s height theorem gives \(\dim\overline B\le k\).

Apply this observation to
$\widehat A=\widehat{\mathcal O}_{V_P,p}$ and to
$\widehat A/J_P\widehat A$.
As before, we identify $\widehat A$ with the completion of the
analytic local ring, so the $u_j$ are elements of $\widehat A$.
If $\widehat{\mathfrak m}$ denotes its maximal ideal, then
\[u'_j-u_j\in\widehat{\mathfrak m}^{\,q}.\]
The original graph satisfies
\[ \dim\widehat A/(u_1,\ldots,u_b)=n,\qquad \dim\widehat A/(J_P\widehat A,u_1,\ldots,u_b)\leq n-2.\]
For all sufficiently large $q$, the observation gives
\[\dim\widehat A/(u'_1,\ldots,u'_b)\leq n,\qquad \dim\widehat A/(J_P\widehat A,u'_1,\ldots,u'_b)\leq n-2.\]
Since $\dim\widehat A=n+b$, quotienting by $b$ elements gives the
opposite inequality for the first quotient. Thus
\[\dim\widehat A/(u'_1,\ldots,u'_b)=n.\]
Completion preserves dimensions. Since
$\mathcal O_{V_P,p}$ is Cohen--Macaulay, the equations
$u'_1,\ldots,u'_b$ form a regular sequence (see \cite[Tag~02JN]{Stacks}).
The germ $(V',p)$ is therefore Gorenstein of dimension $n$. If
$J'=J_P\mathcal O_{V'}$, then
\( \dim\operatorname{Cosupp}(J')\leq n-2.\)

Using properness of $f'$ and shrinking $V'$ equivariantly around $p$,
we may assume that $Y'$ is smooth, that the restricted marked divisors
have simple normal crossings everywhere, and that the preceding
Gorenstein and dimension assertions hold on this neighbourhood.
The morphism $f'$ is an isomorphism over
$V'\setminus\operatorname{Cosupp}(J')$, since this is true of $f_P$
away from $\operatorname{Cosupp}(J_P)$ and is preserved by base change.
Consequently $V'$ is regular in codimension at most one.
Since it is Gorenstein, it is $S_2$ and by Serre's criterion, $V'$ is normal.

We now prove that formation of the blowup of \(J_P\) commutes with base change along the closed embedding \(V'\hookrightarrow V_P\). That is, setting \(J'=J_P\mathcal O_{V'}\), we prove that
\[Y'=W_P\times_{V_P}V'\cong\operatorname{Bl}_{J'}V'.\]
Put
\(\Delta'=\sum_jm_jE'_j.\)
Pulling back the descended generator factorization along \(j:Y'\hookrightarrow W_P\) gives a surjection
\[\mathcal O_{Y'}^{\oplus e}\twoheadrightarrow j^*\mathcal O_{W_P}\!\left(-\sum_jm_jD_j\right)\simeq\mathcal O_{Y'}(-\Delta').\]
Its composite with the canonical inclusion into \(\mathcal O_{Y'}\) is the generator map for \(J'\cdot \mathcal O_{Y'}\). Consequently,
\(J'\cdot \mathcal O_{Y'}=\mathcal O_{Y'}(-\Delta').\)
The natural surjection of graded algebras
\[ \left(\bigoplus_{m\geq0}J_P^m\right) \otimes_{\mathcal O_{V_P}}\mathcal O_{V'} \longrightarrow \bigoplus_{m\geq0}(J')^m\]
induces a closed immersion
\[ \iota:\operatorname{Bl}_{J'}V' \hookrightarrow W_P\times_{V_P}V'=Y'\]
(see \cite[Tag~080E]{Stacks}).
This is an isomorphism over $V'\setminus\operatorname{Cosupp}(J')$.
Its image therefore contains $Y'\setminus\operatorname{Supp}(\Delta')$, which is dense in every
irreducible component of $Y'$.
Since $Y'$ is reduced, the defining ideal of $\iota$ is zero.
Hence \(Y'=\operatorname{Bl}_{J'}V'.\)
In particular, $f'$ is projective and birational. Every component of $\Delta'$ maps into the codimension-at-least-two subset $\operatorname{Cosupp}(J')$, so it is exceptional. 
Conversely, there are no exceptional divisors outside this support.
Thus $f'$ has smooth source and simple normal crossing exceptional
divisor, and is a $G$-equivariant projective log resolution.

It remains to compare discrepancies. Set
\(L=\omega_{W_P}\otimes f_P^*\omega_{V_P}^{-1},\)
and let \(U\subset W_P\) be the open subset over which \(f_P\) is an isomorphism. The natural identification
\(\omega_{W_P}|_U\simeq (f_P^*\omega_{V_P})|_U\)
induces a trivialization \(L|_U\simeq\mathcal O_U\). Let \(\sigma\) be the unique rational section of \(L\) whose restriction to \(U\) corresponds to \(1\) under this trivialization.
Its divisor is supported on the marked divisors $D_j$.
For either graph $Z\in\{V,V'\}$, write
$Y_Z=W_P\times_{V_P}Z$ and $f_Z:Y_Z\to Z$.
Both graph immersions are defined by regular sequences.
Their conormal bundles are identified by pullback, using the same
ordered graph equations on the base and source.
The conormal determinant factors in adjunction therefore cancel,
giving
\[\omega_{Y_Z}\otimes f_Z^*\omega_Z^{-1} \simeq L|_{Y_Z}.\]
On the common isomorphism locus this identification sends
$\sigma|_{Y_Z}$ to the canonical relative section of $f_Z$.

Locally near a point of $F$, choose equations $s_j$ for the marked
divisors through that point and trivialize $L$. We can then write
\(\sigma=\varepsilon\prod_j s_j^{c_j},\)
where $\varepsilon$ is a unit and the $c_j$ are integers.
Transversality ensures that $s_j|_Y$ and $s_j|_{Y'}$ are reduced local
equations for the corresponding exceptional divisors.
It follows that the coefficients of these divisors in
$K_Y-f^*K_V$ and $K_{Y'}-(f')^*K_{V'}$ are both $c_j$.
Since $(V,o)$ is Gorenstein klt, these coefficients are integers
strictly greater than $-1$, and hence are nonnegative.
Thus all exceptional coefficients of $f'$ are nonnegative, and $(V',o')$ is canonical, in particular Gorenstein klt.

Finally, substitution of $t=a(z)$ and $t=a'(z)$ in the descended
equations yields a $G$-equivariant identification
\[
 \mathcal O_{V,o}/\mathfrak m_o^q
 \simeq
 \mathcal O_{V',o'}/\mathfrak m_{o'}^q.
\]
If the action of $G$ on $(V,o)$ is faithful, every nonidentity element
acts nontrivially on some sufficiently high infinitesimal quotient, by Krull intersection.
Since $G$ is finite, one common quotient detects all of its
nonidentity elements.
Choosing $q$ also large enough for this condition makes the induced
action on $(V',o')$ faithful.

\end{proof}

\begin{proposition}\label{p-4.3} Assume Theorem \ref{t-4.1}. Let $(V_i,o_i)$ be a sequence of Gorenstein klt
complex analytic germs of fixed dimension $n\geq2$, and let $g_i$ be
finite-order automorphisms fixing $o_i$. After passing to a subsequence, there are constants $C>0$ and $\ell'\in\mathbb Z_{>0}$, independent of $i$, with the following
properties. For every $i$, there exist a $\langle g_i\rangle$-equivariant projective log resolution
$f_i:Y_i\to V_i$, a point $y_i\in f_i^{-1}(o_i)$ fixed by $h_i:=g_i^{\ell'}$, and a monomial valuation $w_i$ 
with positive rational weights, such that
\[A_{V_i}(w_i)\leq C.\]
The valuation $w_i$ is a positive rational multiple of a prime divisorial valuation centred at $o_i$. For every nonzero meromorphic germ $F$ at
$y_i$ satisfying $h_i^*F=\chi_F(h_i)F$, one has
\[
 w_i(F)\equiv\operatorname{angle}(\chi_F(h_i))\pmod{\mathbb Z}.
\]
In particular, $w_i$ is integer valued on the $ g_i$-invariant field of
meromorphic germs of $(V_i,o_i)$.
\end{proposition}
\begin{proof} 
Put $G_i=\langle g_i\rangle$, regarded as a group of automorphisms of
$(V_i,o_i)$, so that its action is faithful. By Proposition \ref{p-4.2}, we have an algebraic
Gorenstein klt germ $(V_i',o_i')$ with a faithful $G_i$-action and
$G_i$-equivariant projective log resolutions
\[f_i:Y_i\longrightarrow V_i, \qquad f_i':Y_i'\longrightarrow V_i'\]
with the prescribed identifications over $o_i$ and $o_i'$.

We first choose a \(G_i\)-semi-invariant rational section \(\theta_i\) of \(\omega_{V_i'}\) that is regular and nowhere vanishing on a neighbourhood of \(o_i'\).
Indeed, the action on the one-dimensional fibre $\omega_{V_i'}\otimes\mathbb C(o_i')$ is given by a character $\chi_i$. If $\eta_i$ is a local algebraic generator, then
\[\theta_i:=\frac{1}{|G_i|}\sum_{g\in G_i}\chi_i(g)^{-1}g^*\eta_i\] is a non-vanishing semi-invariant section.
Therefore, after shrinking to a $G_i$-invariant
affine neighbourhood of $o_i'$, $\theta_i$ is a generator of \(\omega_{V_i'}\). In particular,
$A_{\theta_i}(w)=A_{V_i'}(w)$ for valuations $w$ centred at $o_i'$.

By Theorem \ref{t-4.1}, after passing to a subsequence, we obtain constants $C_0>0$ and
$\ell\in\mathbb Z_{>0}$, independent of $i$, and valuations
\[w_i'=c_i v_{E_i},\qquad c_i>0,\qquad v_{E_i}=\operatorname{ord}_{E_i}, \]
centred at $o_i'$, where $E_i$ is a prime divisor over $V_i'$ and $v_{E_i}(\mathbb C(V_i')^\times)=\mathbb Z$, such that
\[A_{V_i'}(w_i')\leq C_0\qquad {\rm and}\qquad  w_i'(F)\equiv_\Z \ell\,\operatorname{angle}(\chi_F(g_i))\]
for every nonzero rational $G_i$-eigenfunction $F$ of character $\chi_F$. The restriction of $w_i'$ to the invariant field is
integer valued.

Let $L_i=\mathbb C(V_i')$, let $K_i=L_i^{G_i}$, and let $I_i$ be the inertia subgroup of $v_{E_i}$ for the finite Galois extension $L_i/K_i$. Set $e_i=|I_i|$. Since the residue characteristic is zero, $e_i$ is the ramification index. Thus \[ v_{E_i}(K_i^\times)=e_i\mathbb Z.\]
The integrality of $w_i'|_{K_i}$ implies that \( k_i:=c_i e_i\in\mathbb Z_{>0},\) and in particular, $c_i\in\mathbb Q_{>0}$. 

Choose $\tau_i\in I_i$ whose character on the conormal line
$\mathfrak m_{v_{E_i}}/\mathfrak m_{v_{E_i}}^2$ of the valuation
ring is $\exp(2\pi {\rm i}/e_i)$.
Inertia acts trivially on the residue field $\kappa(E_i)$, so its action on the leading term of a rational $G_i$-eigenfunction $F$ gives
\[ \chi_F(\tau_i) =\exp\!\left(\frac{2\pi {\rm i}\,v_{E_i}(F)}{e_i}\right). \]
The normal basis theorem for the cyclic Galois extension $L_i/K_i$ gives a nonzero rational eigenfunction $F$ whose $G_i$-character is faithful. For this choice, the congruence from Theorem \ref{t-4.1} yields \[\chi_F(g_i^\ell) =\exp(2\pi {\rm i} c_i v_{E_i}(F)) =\chi_F(\tau_i^{k_i}).\]
Since $\chi_F$ is faithful, it follows that \(g_i^\ell=\tau_i^{k_i}\in I_i.\)

Let $Z_i=c_{Y_i'}(E_i)$, so that
$Z_i\subset(f_i')^{-1}(o_i')$.
The inertia group preserves $Z_i$, and the residue-field inclusion
$\mathbb C(Z_i)\hookrightarrow\kappa(E_i)$ shows that it acts trivially on $\mathbb C(Z_i)$, hence on $Z_i$.
Thus $g_i^\ell$ fixes $Z_i$ pointwise.
Choose a general closed point $y_i'\in Z_i$ avoiding every exceptional
component that does not contain $Z_i$. Let $D_{i1}',\ldots,D_{is_i}'$ be the exceptional components through
$y_i'$. They have simple normal crossings, so $s_i\leq n$.
The element $g_i^\ell$ permutes these components. Consequently
\(h_i:=g_i^{\ell n!}=\tau_i^{k_i n!} \) fixes $y_i'$ and preserves every $D_{ij}'$.

First choose local defining equations for \(D'_{i1},\ldots,D'_{is_i}\), and use character averaging to obtain \(h_i\)-eigenfunctions \(z'_{i1},\ldots,z'_{is_i}\) defining these divisors. Choose additional \(h_i\)-eigenvectors so that, together with the differentials of \(z'_{i1},\ldots,z'_{is_i}\) at \(y_i'\), they form a basis of the cotangent space \(T^*_{y_i'}Y_i'\). 
For each additional basis vector, choose a regular-function germ representing it and apply character averaging to obtain an \(h_i\)-eigenfunction with the same cotangent class. The resulting functions \(z'_{i1},\ldots,z'_{in}\) form the required algebraic regular parameters.

Write
\[  h_i^*z_{ij}'=\exp(2\pi {\rm i}\alpha_{ij})z_{ij}', \qquad \alpha_{ij}\in(0,1]\cap\mathbb Q,\]
where the trivial character is represented by $\alpha_{ij}=1$.

On the degree-$m$ piece of the associated graded field of $v_{E_i}$, the element $\tau_i$ acts as multiplication by
$\exp(2\pi {\rm i} m/e_i)$. Hence $h_i=\tau_i^{k_i n!}$ acts on that piece by $\exp(2\pi {\rm i} k_i n!m/e_i)$.
Comparison with the leading term of the $h_i$-eigenfunction $z_{ij}'$ therefore gives
\[ \alpha_{ij} \equiv\frac{k_i n!}{e_i}v_{E_i}(z_{ij}')=n!w_i'(z_{ij}')\pmod{\mathbb Z}.\]
For $1\leq j\leq s_i$, the component $D_{ij}'$ contains $Z_i$, so $w_i'(z_{ij}')>0$. Since $0<\alpha_{ij}\leq1$, we obtain
\[ \alpha_{ij}\leq n!w_i'(z_{ij}') \qquad(1\leq j\leq s_i).\]

Let $y_i\in Y_i$ correspond to $y_i'$, and let $D_{i1},\ldots,D_{is_i}$ be the corresponding exceptional divisors.
The point $y_i$ and these divisors are preserved by $h_i$. Restrict the representation comparisons in Proposition \ref{p-4.2} to
$\langle h_i\rangle$.
Since \(h_i\) fixes \(y_i\) and \(y_i'\) and preserves the corresponding exceptional divisors \(D_{ij}\) and \(D'_{ij}\), Proposition \ref{p-4.2} shows that it acts with the same character on their conormal lines.
The cotangent spaces \(T^*_{y_i}Y_i\) and \(T^*_{y_i'}Y_i'\) are also isomorphic as representations of \(\langle h_i\rangle\). 
Because the marked divisors have simple normal crossings, their conormal lines form a direct sum in each cotangent space. 
By semisimplicity, these sums admit invariant complements, and the two complementary representations have the same characters, counted with multiplicity.
We can therefore choose analytic $h_i$-eigenparameters $z_{i1},\ldots,z_{in}$ at $y_i$ with
\[ D_{ij}=(z_{ij}=0)\quad(1\leq j\leq s_i),\qquad h_i^*z_{ij}=\exp(2\pi {\rm i}\alpha_{ij})z_{ij}.\]
Let $w_i$ be the monomial valuation with weights $w_i(z_{ij})=\alpha_{ij}$.
These weights are positive and rational, so $w_i$ is a positive rational multiple of a prime-divisor order, realized by the exceptional divisor of the corresponding weighted blowup at $y_i$. Its centre on $V_i$ is $o_i$.

Let $a_{ij}$ be the coefficient of $D_{ij}$ in $K_{Y_i}-f_i^*K_{V_i}$. By Proposition \ref{p-4.2}, it is also the coefficient of $D_{ij}'$ in $K_{Y_i'}-(f_i')^*K_{V_i'}$.
Since the germs are Gorenstein klt, these coefficients are integers greater than $-1$, hence $a_{ij}\geq0$.
Our choice of $y_i'$ ensures that $D_{i1}',\ldots,D_{is_i}'$ are precisely the exceptional components containing $Z_i$.
The discrepancy transformation formula therefore gives \[A_{V_i'}(w_i') =A_{Y_i'}(w_i')+   \sum_{j=1}^{s_i}a_{ij}w_i'(z_{ij}').\]
As $Y_i'$ is smooth, $A_{Y_i'}(w_i')>0$.
Using the discrepancy formula for a monomial valuation on $Y_i$, we conclude that
\[  A_{V_i}(w_i) =\sum_{j=1}^n\alpha_{ij}   +\sum_{j=1}^{s_i}a_{ij}\alpha_{ij}\leq n+n!\sum_{j=1}^{s_i}a_{ij}w_i'(z_{ij}')\leq n+n!A_{V_i'}(w_i')\leq n+n!C_0.  \]

It remains to prove the character congruence for analytic
meromorphic germs at $y_i$.
If $F$ is a nonzero holomorphic $h_i$-eigenfunction, then every nonzero
monomial in its convergent expansion in the $z_{ij}$ has character
$\chi_F(h_i)$. The weight of each such monomial is therefore
congruent to $\operatorname{angle}(\chi_F(h_i))$ modulo $\mathbb Z$.
The positive rational weights ensure that their minimum is attained,
so the same congruence holds for $w_i(F)$.

For a nonzero meromorphic eigenfunction $F=a/b$, with $a,b$
holomorphic and $b\neq0$, let
\( B:=\prod_{\tau\in\langle h_i\rangle}\tau^*b.\)
Then $B$ is a nonzero invariant holomorphic germ, and $FB$ is
holomorphic with the same character as $F$.
Subtracting the congruences for $FB$ and $B$ proves the desired congruence for $F$.
Finally, a $G_i$-invariant meromorphic germ on $V_i$ pulls back to
an $h_i$-invariant germ at $y_i$, so its value under $w_i$ is an
integer. The required constants are
\[ \ell'=\ell n!,\qquad C=n+n!C_0, \]
which are independent of $i$ on the given subsequence.

\end{proof} 
\section{Cartier index bounds}

We will need the following well known lemma.
\begin{lemma}\label{lem:gorenstein-index-cover}
Let $(X,B)$ be a strongly $\Q$-factorial klt germ, $dB$ be integral for some positive integer $d$, and set
\(L=d(K_X+B).\)
Let $r$ be the Cartier index of the integral $\mathbb Q$-Cartier
divisor $L$ at the closed point $x$. Then there exists a finite
quasi-\'etale Galois morphism of pointed germs
\[\pi:(V,o)\longrightarrow (X,x)\]
with group $H$ such that
\begin{enumerate}
\item $H$ fixes $o$,
\item $V$ is Gorenstein klt,
\item $K_V$ and $\pi^{*}L$ are Cartier,
\item $\pi^{*}L$ admits an $H$-semi-invariant local generator
      $\eta$, and its character $\chi:H\to\mathbb C^*$ has image
      $\mu_r$.
\end{enumerate}
\end{lemma}
\begin{proof} The normalized local index-$r$ cover of $L$ is cyclic,
quasi-\'etale, and its tautological generator has primitive
$\mu_r$-character.  We then take a canonical index-one cover of this
germ and the corresponding Galois closure over $(X,x)$. Let $\widetilde V\to X$ be the finite quasi-\'etale Galois cover
constructed above, with Galois group $G$. 
Choose \(o\in\widetilde V\) over \(x\), set \(H:=\operatorname{Stab}_G(o)\), and put \((V,o):=(\widetilde V,o)\). Then
\( \pi:(V,o)\longrightarrow(X,x) \)
is a finite quasi-étale \(H\)-Galois cover of analytic germs. In particular, \(H\) fixes \(o\), and every meromorphic germ pulled back from \((X,x)\) is \(H\)-invariant.  Since the original index-$r$ cover has a
unique point over $x$, the map $H\to\mu_r$ is surjective.  Hence the pullback of the tautological generator has character of exact order $r$.
All covers are quasi-\'etale.  Therefore $V$ is klt and \(K_V=\pi^*K_X\)
is Cartier.  Thus $V$ is Gorenstein.  Moreover,
\(\pi^{*}(dB)=\pi^{*}L-dK_V \) is Cartier.
\end{proof}
\subsection{Index covers and discrepancy residues}\label{ss-5.1} We fix a positive integer $d$ such that $dB$ is a $\Z$-divisor and let $L_i=dD_i=d(K_{X_i}+B_i)$ so that each $L_i$ is an integral
$\Q$-Cartier divisor. 
Let $r_i$ be the Cartier index of $L_i$ at a point $x_i\in X_i$. Let $\pi_i:(V_i,o_i)\to (X_i,x_i)$ be the germ of the Galois cover given by Lemma \ref{lem:gorenstein-index-cover}. Thus $V_i\to X_i$ is a finite quasi-\'etale, Galois cover with Galois group $H_i$, $V_i$ is normal $\pi_i ^*L_i$ is Cartier, $K_{V_i}$ is Cartier, and \(\eta_i\) denotes the semi-invariant generator of $\pi_i ^{*}L_i$. We write $K_{V_i}+B_{V_i}=\pi_i ^*(K_{X_i}+B_i)$, then $(V_i,B_{V_i})$ is klt and $V_i$ is Gorenstein. Therefore $V_i$ has canonical singularities. 
Note that $\OO_{V_i}(d(K_{V_i}+B_{V_i}))=\OO_{V_i}(\pi_i ^* L_i)$ and hence $\OO_{V_i}(d(K_{V_i}+B_{V_i}))$ is $H_i$-linearized. Thus \(H_i\) acts on the one dimensional fibre at \(o_i\)
via a character $\chi_i:H_i\to \C ^*$. We may assume that $h^*(\eta _i)=\chi_i (h)\cdot \eta _i$ for $h\in H_i$. 

Fix an arbitrary sequence of elements \(g_i\in H_i\). Applying Proposition \ref{p-4.3} to the sequence \((V_i,o_i,g_i)\), and passing to a subsequence, we obtain valuations \(w_i\), a constant \(C>0\), and an integer \(\ell'>0\), with \(C\) and \(\ell'\) independent of \(i\). The particular choice of \(g_i\) will be made in the proof of Proposition \ref{p-5.1}.
Let $v_i$ be the restriction of $w_i$ to the field of meromorphic functions on $(X_i,x_i)$. Since the pull-back of these functions are $H_i$ invariant, then $v_i$ is integer valued. 
Recall that $A_i(\ldots ):=A_{X_i,B_i}(\ldots )$. We have \begin{equation}\label{e-9}A_i(v_i)=A_{V_i,B_{V_i}}(w_i)\leq A_{V_i}(w_i).\end{equation}
Let $y_i$ be the point on the equivariant projective
log resolution of $V_i$ given by Proposition \ref{p-4.3}, and $z_1, \ldots , z_n$ the eigen-coordinates centered at $y_i$. Let $F(dz_1\wedge \ldots \wedge dz_n)^{\otimes d}$ be the pullback of $\eta_i$.
The divisor corresponding to $\eta _i$ is $-dB_{V_i}$. The function $F$ is semi-invariant under the element
$g_i^{\ell '}$. The monomial discrepancy formula gives \[dA_i(v_i)=w_i(F)+d\sum _{j=1}^nw_i(z_j).\] 
Writing character angles with $h_i=g_i^{\ell'}$, we have $h_i^*(f)=\chi _f(h_i)\cdot f$ and by Proposition \ref{p-4.3}, $w _i(f)\equiv _\Z {\rm arg}\ \chi _f(h_i)$. Thus we obtain
\begin{equation}\label{e-10}dA_i(v_i)\equiv _\Z  {\rm arg}\ \chi _F(h_i)+d \sum _{j=1}^n{\rm arg}\ \chi _{z_j}(h_i) \equiv _\Z {\rm arg} \left(\chi _F\cdot (\chi _{z_1}\cdots \chi _{z_n})^d\right)(h_i)=\end{equation} \begin{equation*}= {\rm arg}\  \chi_i(h_i) = {\rm arg}\  \chi_i(g_i^{\ell'}) \equiv _\Z \ell'{\rm arg}\ \chi_i(g_i).\end{equation*}
\subsection{Realization on a compact modification}\label{s-5.2} We wish to apply Proposition \ref{lem:countable-limits} to our situation. However since the above construction $(V_i,o_i)\to (X_i,x_i)$ is local, and Proposition \ref{lem:countable-limits} applies to compact varieties, we must compactify our setting.

Write \(A\) for the convergent local ring of \(X_i\) at \(x_i\), $V\to {\rm Specan}(A)$ the corresponding cover
and $v,w$ for $v_i,w_i$. Let \(K=\operatorname{Frac}(A)\), and let \(\mathcal O_v\subset K\) be the valuation ring of \(v\), with maximal ideal \(\mathfrak m_v\) and residue field \(\kappa(v)=\mathcal O_v/\mathfrak m_v\). The residue fields of $w$ and $v$ have transcendence degree $n - 1$ over $\C$. Since \(K(V)/K\) is finite, the residue fields of \(w\) and of its restriction \(v\) have the same transcendence degree \(n-1\) over \(\mathbb C\).

Choose $f,h_1, \ldots , h_{n-1}\in A$ such that $v(f)=v(h_1)=\ldots =v(h_{n-1})>0$, and  the residue classes of 
\(h_1/f,\ldots,h_{n-1}/f\)
 in the residue field \(\kappa(v)\) are algebraically independent over $\C$.
 Fix $N>0$ such that $v(\mathfrak m^N)>v(f)$ and let
\[J=(f,h_1, \ldots , h_{n-1})+\mathfrak m^N, \qquad T={\rm Bl}_J({\rm Specan}(A)).\]
Since $v(\mathfrak m^N)>v(f)$, the image of $f^{-1}(J)\cap \OO _v$ in $\kappa (v)$ contains the algebraically independent elements $h_1/f,\ldots,h_{n-1}/f$ and hence the center $C$ of $v$ on $T$ has dimension $n-1$. 
On the normalization \(T^\nu\), the centre of \(v\) is therefore a prime divisor \(E\), and 
\(v=q\,\operatorname{ord}_E,\) \(q\in\mathbb Z_{>0}.\)
The ideal  \(J\) is \(\mathfrak m\)-primary and glueing it with $\mathcal O_{X_i}$ away from $x_i$ we obtain a coherent ideal sheaf \(\mathfrak a\) co-supported at \(x_i\).  The normalization of the blow up of \(\mathfrak a\) defines a morphism $Y_i\to X_i$ which is an isomorphism over $X_i\setminus x_i$ and whose germ at $x_i$ coincides with $T^\nu\to {\rm Specan}(A)$.
Since the $X_i$ belong to a given mmp, they are all birational and hence we may consider the corresponding valuations $v_i$ centered at $x_i$ as valuations corresponding to divisors over the fixed resolution \(W\) of $X$.
\subsection{Boundedness along the program}
\begin{proposition}\label{p-5.1} 
Let \((X,B)\) be a strongly \(\mathbb Q\)-factorial compact Kähler klt pair, where \(B\) is an effective \(\mathbb Q\)-divisor. Assume that \(K_X+B\) is pseudo-effective, and let \(\omega\) be a Kähler class such that \(K_X+B+\omega\) is nef. Fix a sequence
\[
(X,B)=(X_0,B_0)\dashrightarrow(X_1,B_1)\dashrightarrow
(X_2,B_2)\dashrightarrow\cdots
\]
of steps of the \((K_X+B)\)-minimal model program with scaling of \(\omega\), where \(B_i\) is the pushforward of \(B\), and set
\(D_i:=K_{X_i}+B_i. \)
Then there exists an integer \(N>0\), independent of \(i\), such that \(ND_i\) is Cartier on \(X_i\) for every \(i\geq0\).
\end{proposition}
\begin{proof} We may assume that $n=\dim X\geq 2$ and the mmp $X\dasharrow X_1\dasharrow X_2\dasharrow \ldots $ is infinite. Fix \(d\) as in Subsection \ref{ss-5.1}, and put \(L_i=dD_i\). Suppose that the local Cartier indices of the \(L_i\) are unbounded.  We can therefore choose an infinite subsequence of indices, and points \(x_i\in X_i\), such that
\(r_i:=\operatorname{ind}_{x_i}(L_i)\) and \(\lim r_i=\infty.\)
All subsequent indices \(i\) refer to this subsequence of the original minimal model program. Let \((V_i,o_i)\to(X_i,x_i)\) and \(H_i\) be the covers and groups of Lemma \ref{lem:gorenstein-index-cover}. 

The generator characters \(\chi_i:H_i\longrightarrow\mu_{r_i} \) are surjective.
Fix \(\alpha\in[0,1)\), set
\(k_i:=\lfloor r_i\alpha\rfloor \)
and choose
\[g_i\in\chi_i^{-1}
\left(\exp\left(2\pi i\frac{k_i}{r_i}\right)\right),\]
which is possible because \(\chi_i:H_i\to\mu_{r_i}\) is surjective. Then
\[ \arg\chi_i(g_i)=\frac{k_i}{r_i}, \qquad 0\leq\alpha-\frac{k_i}{r_i}<\frac1{r_i},\]
and so \(\arg\chi_i(g_i)\to\alpha\).

By Proposition \ref{p-4.3}, passing to a subsequence, there exists $C>0$ and $\ell'\in \Z^{>0}$ such that,
after passing to a further subsequence, we have $A_{V_i}(w_i)\leq C$ and for every nonzero meromorphic germ $f$ at $y_i$ with $g^{\ell'}_i$-semi-invariant character $\chi_f$, then
$w_i(f ) \equiv _\Z {\rm arg}\ \chi_f (g^{\ell'}_i )$.
By \eqref{e-8}, $A_i(v_i)\geq A_0(v_i)\geq \epsilon A_W(v_i)$ (here $\epsilon >0$ is determined by the given mmp).  By \eqref{e-9}, $A_i(v_i)\leq A_{V_i}(w_i)$. Putting this together we have\[ \epsilon A_{W}(v_i)\leq A_i(v_i)\leq A_{V_i}(w_i)\leq C\] so that $A_W(v_i)\leq \frac C \epsilon$
is uniformly bounded. Since \(E_i(v_i)=A_i(v_i)-A_0(v_i)\geq 0,\) and since $(X,B)$ is klt, we have $A_0(v_i)>0$, then \[0\leq E_i(v_i)\leq A_i(v_i)\leq C\]
and passing to a subsequence we may assume that $E_i(v_i)$ converges to \[\lim _{i\to \infty}E_i(v_i)=e\in \mathcal E\] where $\mathcal E$ is the countable set given in Proposition \ref{lem:countable-limits}.

We now choose a positive integer $M_0$
divisible by $d$ such that $M_0D_0$ is Cartier. Since $v_i$ is integer-valued,
$M_0A_0(v_i)$ is an integer and so $M_0E_i(v_i)\equiv _\Z M_0A_i(v_i)$. 
For any real number $r$, we denote $[r] \in\R/\Z$. 
Multiplying \eqref{e-10} 
by $M_0/d$ and passing to the limit of the corresponding residues in $\R/\Z$, we obtain
\[ [M_0e]=\lim_{i\to \infty} [M_0E_i(v_i)] =\lim_{i\to \infty} [M_0A_i(v_i)]= \lim_{i\to \infty} [\frac{M_0}d \ell'{\rm arg}\ \chi_i(g_i)]=[\frac{M_0}d \ell'\alpha] .\]

The set $\mathcal E$ and the integers $M_0$, $d$ are independent of the choice of $\alpha$. The exponent $\ell'$ may depend on the choice of $\alpha$, but is
always a positive integer. For any fixed $e \in \mathcal E$ and $\ell'\in \Z_{>0}$, the congruence
allows only finitely many $\alpha$'s. Taking the union over the countable set $\mathcal E$ and
all positive integers $\ell'$ gives a countable set. Our choice of $\alpha$ was arbitrary
in the uncountable set $\R/\Z$. This is impossible and so the 
local indices $r_i$ of the $L_i$ are bounded. If $r$ is their least common multiple and $N=dr$, then $ND_i$ is Cartier for all $i\geq 0$.
\end{proof}

\section{Proof of the main result}
\begin{proof}[Proof of Theorem \ref{t-main}] We denote by $X\dasharrow X_1\dasharrow X_2 \dasharrow \ldots $ a sequence of steps of the $K_X+B$ mmp with scaling of $\omega$, $\phi _n:X\dasharrow X_n$ and $B_n=(\phi_n) _*B$. As discussed in Subsection \ref{ss-sketch}, we may assume that $K_X+B$ is pseudo-effective, $\lambda _i>0$ and $\lim \lambda _i=0$. We may assume that there are only finitely many divisorial contractions and hence,  that each $\psi _i:X_i\dasharrow X_{i+1}$ is a flip for $i\gg 0$. Let $f_i:X_i\to Z_i$ be the corresponding flipping  contraction and $f_i^+:X_{i+1}\to Z_i$ the corresponding flipped contraction. We also let $p_i:U_i\to X_i$, $q_i:U_i\to X_{i+1}$ be a common projective log resolution of $(X_i,B_i)$ and $(X_{i+1},B_{i+1})$, \[S_p:=K_{U_i}-{p_i^*}(K_{X_i}+B_i), \qquad S_q=K_{U_i}-{q_i^*}(K_{X_{i+1}}+B_{i+1})\] \[\Phi_i=S_q-S_p=p_i^*(K_{X_{i}}+B_{i})-q_i^*(K_{X_{i+1}}+B_{i+1}).\] The $\Q$-divisor $\Phi_i$ is  effective and $p_i,q_i$-exceptional.
Note that by Proposition \ref{p-5.1} there exists an integer $R\geq 2$ such that each $m(K_{X_i}+B_i)$ is Cartier for any positive integer $m$ divisible by $R$. Note that it follows that $m\Phi_i$ is also Cartier. For $m$ divisible by $R$, we let
\[\LL_{m,i}=q_i^*\OO_{X_{i+1}}(m(K_{X_{i+1}}+B_{i+1}))\otimes \OO _{U_i}(\lceil S_q\rceil).\] Note that \[q_i^*(m(K_{X_{i+1}}+B_{i+1}))+\lceil S_q\rceil\equiv K_{U_i}+\{-S_q\}+(m-1)q_i^*(K_{X_{i+1}}+B_{i+1})\] where $(U_i,\{-S_q\})$ is klt and $(m-1)q_i^*(K_{X_{i+1}}+B_{i+1})$ is nef and big over $Z_i$ (as $K_{X_{i+1}}+B_{i+1}$ is ample over $Z_i$ and $U_i\to Z_i$ is birational). Let $g_i=f_i\circ p_i=f_{i}^+\circ q_i$. By relative Kawamata–Viehweg vanishing (see \cite[Lemma 2.13, Remark 2.14]{Fuj23}) we have \[R^j{p_i}_*\LL_{m,i}=0,\qquad R^j{g_i}_*\LL_{m,i}=0,\qquad \forall j>0.\]
Since $\lceil S_q\rceil-m\Phi_i=\lceil S_q-m\Phi_i\rceil=\lceil S_p-(m-1)\Phi_i\rceil$, then\[q_i^*(m(K_{X_{i+1}}+B_{i+1}))+\lceil S_q\rceil\sim p_i^*(m(K_{X_{i}}+B_{i}))+\lceil S_p-(m-1)\Phi_i\rceil.\]
\begin{claim} Fix $m$ a positive integer divisible by $R$. The divisor  \[\lceil S_p-(m-1)\Phi_i\rceil\] is effective and $p_i$-exceptional for all $i\gg 0$.
\end{claim}
\begin{proof} If $v$ is a divisor on $X_i$, then $\Phi_i(v)=0$ and $S_p(v)=-{\rm mult}_v(B_i)>-1$ and hence ${\rm mult}_v(\lceil S_p-(m-1)\Phi_i\rceil)=0$, i.e. $\lceil S_p-(m-1)\Phi_i\rceil$ is $p_i$-exceptional.

It suffices to show  that
\[A_i(v)-(m-1)\Delta _i(v)={\rm mult}_{E_v}(S_p-(m-1)\Phi_i)+1>0\]
for any $p_i$-exceptional divisorial valuation $v$.

By \eqref{e-8}, there is an $\epsilon >0$ such that $A_i(v)\geq \epsilon A_W(v)$. Fix $s>(m-1)/\epsilon$. By  Proposition \ref{p-5} we may assume that for all $i\gg 0$, we have $\Delta _i(v)\leq A_W(v)/s$ for any divisorial valuation $v$. Therefore,   \[A_i(v)-(m-1)\Delta _i(v)\geq (1-\frac{m-1}{\epsilon s})A_i(v)>0.\]
For $v={\rm ord}_E$ this implies that $1+{\rm mult}_E(S_{p}-(m-1)\Phi _i)>0$ and hence every coefficient of $ S_p-(m-1)\Phi_i $ is greater than $-1$ so that $\lceil S_p-(m-1)\Phi_i\rceil$ is effective.
\end{proof}
 Since \(\lceil S_p-(m-1)\Phi_i\rceil\) is effective and $p_i$-exceptional ${p_i}_*\LL_{m,i}=\OO _{X_i}(m(K_{X_i}+B_i))$. By the Leray spectral sequence, then \[R^j{f_i}_*\OO _{X_i}(m(K_{X_i}+B_i))=R^j{g_i}_*\mathcal L _{m,i}=0, \qquad {\rm for}\ j>0.\]
\begin{claim}\label{c-1} Let $\mathcal F _i:=\OO _{X_i}(R(K_{X_i}+B_i))$, then $\mathcal F _i$ is $f_i$-generated for $i\gg 0$.
\end{claim}
Grant the claim for the time being. If $\Sigma $ is a flipping curve, then $(K_{X_i}+B_i)\cdot \Sigma\geq 0$ which is a contradiction. Therefore the theorem is proved.
It remains to prove the claim.
\begin{proof}[Proof of Claim \ref{c-1}]
By \cite[Theorem 6.4]{Fuj22} there is an integer $N>0$ depending only on $n:=\dim X$ such that $\mathcal A _i:=\OO _{X_i}(-NR(K_{X_i}+B_i))$ is $f_i$-generated and $f_i$-ample. Let $m_j=R+NRj$ for $1\leq j\leq n$. 
By the above arguments, we may assume that \[R^j{f_i}_*(\mathcal F_i\otimes \mathcal A_i^{-j})=R^j{f_i}_*\OO _{X_i}(m_j(K_{X_i}+B_i))=0,\qquad 1\leq j\leq n.\]
By relative Castelnuovo–Mumford regularity $\mathcal F_i$ is $f_i$-generated (see \cite[Example 1.8.24]{Laz04}).
\end{proof}
\end{proof}
\begin{proof}[Proof of Corollary \ref{c-1.2}] Replacing $(X,B)$ by a small strong $\Q$-factorialization (see \cite[Theorem 3.2]{HLX26}), we may assume that $X$ is strongly $\Q$-factorial. Fix $\omega $ a K\"ahler form such that $K_X+B+\omega$ is K\"ahler. By Theorem \ref{t-main}, after running a mmp with scaling we obtain a minimal model. Therefore, we may assume that $K_X+B$ is nef and it remains to show that $K_X+B$ is semiample.
This follows from \cite[Theorem 1.2]{OAI3}. Note that \cite[Assumption 1.1]{OAI3} holds by \cite[Theorem 1.2]{OAI4} or can be deduced from \cite[Theorem 1.1]{HPS18} (see the discussion following \cite[Theorem 1.2]{OAI4}).
\end{proof}

\end{document}